\documentclass{article}

\usepackage[margin=1in]{geometry}
\usepackage{amsmath, amssymb, tikz-cd, centernot, mathtools, csquotes,
            mathrsfs, hyperref, amsthm, xparse, microtype, enumitem, bbm}
\usepackage[utf8]{inputenc}
\usepackage[backend=biber, style=numeric, backref=true]{biblatex}
\usepackage[nameinlink,capitalize]{cleveref}
\usepackage{IEEEtrantools}

\hypersetup{
  colorlinks,
  linktoc=section,
  linkcolor=blue,
  citecolor=blue,
  urlcolor=blue,
}

\newcounter{smarttheorem:current} 
\newcounter{smarttheorem:next} 

\NewDocumentCommand{\newsmarttheorem}{m o m o}{
  \IfNoValueTF{#2}{\newcounter{#1}}{}
  \newtheorem{hidden#1}[\IfValueTF{#2}{#2}{#1}]{#3}

  \NewDocumentEnvironment{#1}{o}{
    \IfNoValueTF{##1}{\begin{hidden#1}}{\begin{hidden#1}[{##1}]}
      \label{smarttheorem:\arabic{smarttheorem:next}}
      \edef\currentlabel{\arabic{smarttheorem:next}}
      \stepcounter{smarttheorem:next}
  }{
      \setcounter{smarttheorem:current}{\currentlabel}
    \end{hidden#1}
  }
  \Crefname{hidden#1}{#3}{\IfNoValueTF{#4}{#3s}{#4}}
}

\let \vqedsymbol \qedsymbol

\NewDocumentEnvironment{pf}{o o}{
  \IfValueTF{#2}{
    \edef \qedsymbollabel {#2}
    \renewcommand \qedsymbol {\vqedsymbol \, \qedsymbollabel}
  }{
    \renewcommand \qedsymbol \vqedsymbol
  }
  \IfValueTF{#1}{\begin{proof}[#1]}{\begin{proof}}
}{
  \end{proof}
}

\NewDocumentEnvironment{rpf}{O{\proofname}
O{smarttheorem:\arabic{smarttheorem:current}}}{
  \begin{pf}[#1][\noexpand\Cref{#2}]
}{
  \end{pf}
}

\NewDocumentEnvironment{lrpf}{m}{
  \begin{rpf}[Proof of \cref{#1}][#1]
}{
  \end{rpf}
}

\newcommand{\plr}{\item[\((\implies)\)]}
\newcommand{\prl}{\item[\((\impliedby)\)]}

\NewDocumentEnvironment{ea}{O{rCl}}{
  \begin{IEEEeqnarray*}{#1}
}{
  \end{IEEEeqnarray*}
  \ignorespacesafterend
}

\NewDocumentEnvironment{tcd}{s}{
  \IfBooleanTF{#1}{\begin{equation}}{\begin{equation*}}
    \begin{tikzcd}
}{
    \end{tikzcd}
  \IfBooleanTF{#1}{\end{equation}}{\end{equation*}}
  \ignorespacesafterend
}

\newcounter{case}
\newenvironment{caselist}{
  \setcounter{case}{0}\begin{description}
}{
  \end{description}
}

\newcommand \case {\setcounter{case}{\value{case}+1}\item[Case \thecase.]}

\numberwithin{theoremcounter}{section}
\theoremstyle{plain}
\newsmarttheorem{theorem}[theoremcounter]{Theorem}
\newsmarttheorem{lemma}[theoremcounter]{Lemma}
\newsmarttheorem{proposition}[theoremcounter]{Proposition}
\newsmarttheorem{claim}[theoremcounter]{Claim}
\newsmarttheorem{subclaim}[theoremcounter]{Subclaim}
\newsmarttheorem{corollary}[theoremcounter]{Corollary}[Corollaries]
\newsmarttheorem{fact}[theoremcounter]{Fact}
\newsmarttheorem{todo}{TODO}

\theoremstyle{definition}
\newsmarttheorem{definition}[theoremcounter]{Definition}
\newsmarttheorem{notation}[theoremcounter]{Notation}

\theoremstyle{remark}
\newsmarttheorem{question}[theoremcounter]{Question}
\newsmarttheorem{remark}[theoremcounter]{Remark}
\newsmarttheorem{scratch}[theoremcounter]{Rough work}
\newsmarttheorem{aside}[theoremcounter]{Aside}
\newsmarttheorem{conjecture}[theoremcounter]{Conjecture}
\newsmarttheorem{exercise}[theoremcounter]{Exercise}
\newsmarttheorem{example}[theoremcounter]{Example}
\newsmarttheorem{ednote}[theoremcounter]{Editor's note}

\renewcommand \epsilon \varepsilon
\renewcommand \phi \varphi

\renewcommand \th {^\mathrm{th}}
\newcommand \telse {\text{else}}
\newcommand \tif {\text{if }}

\newcommand \lit [1] {\item[#1]}

\def\Indep#1#2{#1\setbox0=\hbox{$#1x$}\kern\wd0\hbox to 0pt{\hss$#1\mid$\hss}
\lower.9\ht0\hbox to 0pt{\hss$#1\smile$\hss}\kern\wd0}

\def\notindep#1#2{#1\setbox0=\hbox{$#1x$}\kern\wd0
\hbox to 0pt{\mathchardef\nn=12854\hss$#1\nn$\kern1.4\wd0\hss}
\hbox to
0pt{\hss$#1\mid$\hss}\lower.9\ht0 \hbox to 0pt{\hss$#1\smile$\hss}\kern\wd0}

\newcommand \tfae {the following are equivalent}
\newcommand \Tfae {The following are equivalent}

\newenvironment{menumerate}{\mbox{}\begin{enumerate}}{\end{enumerate}}

\newenvironment{mdescription}{\mbox{}\begin{description}}{\end{description}}

\DeclarePairedDelimiter \abs \lvert \rvert
\DeclarePairedDelimiter \nrm \lVert \rVert
\DeclarePairedDelimiter \pars ( )

\DeclarePairedDelimiterX \set [1] \{ \} {\, #1 \,}

\newcommand \col [2] {\begin{pmatrix}#1\\#2\end{pmatrix}}
\newcommand \tcol [3] {\begin{pmatrix}#1\\#2\\#3\end{pmatrix}}
\newcommand \fcol [4] {\begin{pmatrix}#1\\#2\\#3\\#4\end{pmatrix}}

\newcommand \tpowmult  {\times \mathord\restriction t^\mathbb{N}}

\DeclareMathOperator \lspan {span}

\DeclareMathOperator \id {id}

\DeclareMathOperator \ndm {End}

\DeclareMathOperator \aut {Aut}
\DeclareMathOperator \thy {Th}

\DeclareMathOperator \lcm {lcm}

\title{Contributions to the theory of \(F\)-automatic sets}
\author{Christopher Hawthorne\thanks{%
This work was partially supported by an NSERC PGS-D and an NSERC CGS-D.}}

\begin{document}
\maketitle
\begin{abstract}
  Fix an abelian group \(\Gamma \) and an injective endomorphism \(F\colon
  \Gamma \to\Gamma \). Improving on the results of \cite{bell19}, new
  characterizations are here obtained for the existence of spanning sets,
  \(F\)-automaticity, and \(F\)-sparsity. The model theoretic status of these
  sets is also investigated, culminating with a combinatorial description of
  the \(F\)-sparse sets that are stable in \((\Gamma ,+)\), and a proof that
  the expansion of \((\Gamma ,+)\) by any \(F\)-sparse set is NIP. These
  methods are also used to show for prime \(p\ge 7\) that the expansion of
  \((\mathbb{F}_p[t],+)\) by multiplication restricted to \(t^\mathbb{N}\) is
  NIP.
\end{abstract}
\tableofcontents
\section{Introduction}
A set \(S\) of natural numbers is said to be \(p\)-automatic if the set of
representations base \(p\) of the elements of \(S\) forms a regular
language---meaning that it is recognized by a finite automaton. Motivated by
the isotrivial Mordell-Lang problem, Bell and Moosa were led in \cite{bell19}
to extend this notion to the context where the natural numbers are replaced by
an abelian group \(\Gamma \) and multiplication by \(p\) is replaced by a fixed
injective endomorphism \(F\) of \(\Gamma \). The intended examples were when
\(\Gamma \) is a finitely generated subgroup of a semiabelian variety \(G\)
over a finite field \(\mathbb{F}_q\), with \(\Gamma \) preserved by the
\(q\)-power Frobenius endomorphism \(F\) of \(G\). In that setting, Bell and
Moosa show in \cite{bell19} that if \(X\) is a closed subvariety of \(G\)
then \(X\cap\Gamma \) is \(F\)-automatic; they do so using the isotrivial
Mordell-Lang theorem of Moosa and Scanlon in \cite{moosa04}, which describes
\(X\cap\Gamma \).  These results have recently been made effective and
generalized to arbitrary commutative algebraic groups by Bell, Ghioca, and
Moosa in \cite{bell20}.

Our interest here is in the general theory of \(F\)-automata. The precise
definition of an \(F\)-automatic subset of \(\Gamma \) is recalled in
\cref{sec:prelims}, but let us give an informal explanation now. First one
needs to know that \((\Gamma ,F)\) admits a ``spanning set'': this is a finite
subset \(\Sigma \) such that every element of \(\Gamma \) can be expressed as
\(s_0+Fs_1+\cdots + F^ns_n\) for some \(s_0,\ldots ,s_n\in\Sigma \). Then a
subset \(S\) of \(\Gamma \) is said to be \(F\)-automatic if the set of words
\(s_0s_1\cdots s_n\) such that \(s_0+Fs_1+\cdots +F^ns_n\in S\) forms a regular
language on the alphabet \(\Sigma \).

Our goals in this paper are twofold: (1) to clarify and further develop the
foundations of this theory, and (2) to investigate the model-theoretic
properties of \(F\)-automatic sets, in particular with respect to the
stability-theoretic hierarchy.

Regarding foundations, our main results are:
\begin{itemize}
  \item
    A characterization of when \((\Gamma ,F)\) admits a spanning set in terms
    of the existence of a height function on \(\Gamma \)
    (\cref{thm:spanlenheight}) and, in the finitely generated case, in terms of
    the eigenvalues of \(F \otimes _{\mathbb{Z}}\id_\mathbb{C}\)
    (\cref{thm:spaneigen}).
  \item
    A characterization of \(F\)-automaticity in terms of kernels; this is
    \cref{thm:cosetdefauto} below. This result brings the basic theory of
    \(F\)-automatic sets closer to the classical case of \(p\)-automaticity.
  \item
    A characterization of \(F\)-sparsity. In the theory of regular languages
    there is a natural sparse/non-sparse dichotomy in terms of the growth of
    the number of accepted words of bounded length. This was extended and
    investigated in \cite{bell19}; see \cref{def:fsparse} below for details. We
    clarify some properties of \(F\)-sparse sets (answering a question in
    \cite{bell19} along the way) and characterize them in terms of length
    functions (\cref{thm:fsparsechar}). As a consequence, we obtain in
    \cref{cor:fsparsesigsparse} a useful criterion for verifying
    \(F\)-sparsity; this criterion is used by Bell-Ghioca-Moosa in the recent
    preprint \cite{bell20} on effective isotrivial Mordell-Lang.
\end{itemize}
We next turn to the model-theoretic analysis. Here our main results are:
\begin{itemize}
  \item
    An explicit characterization of stable \(F\)-sparse sets as finite Boolean
    combinations of translates of finite sums of sets of the form
    \(\set{a+Fa+\cdots +F^na : n <\omega }\)---these are the ``groupless
    \(F\)-sets'' of \cite{moosa04}. That the groupless \(F\)-sets are stable is
    from \cite{moosa04}; the converse is \cref{thm:stablesparsechar} below.
  \item
    The production of some new NIP expansions of \((\Gamma ,+)\); see
    \cref{thm:edp1dim}. This includes \((\Gamma ,+,A)\) for any \(F\)-sparse
    \(A\subseteq \Gamma \), but it also includes some examples that are not
    even \(F\)-automatic, most notably the expansion of \((\mathbb{F}_p[t],+)\)
    by the graph of multiplication restricted to \(t^\mathbb{N}\), when \(p\ge
    7\). This last example is \cref{thm:polysnip} below.
\end{itemize}
My thanks to Luke Franceschini for helping me work through the details of
\cref{eg:luke}. Warm thanks to my advisor, Rahim Moosa, for his excellent
guidance, thorough editing, and many helpful discussions.
\section{Preliminaries}
\label{sec:prelims}
We first recall the relevant theory of regular languages; the interested reader
is directed to \cite{yu97} for further details.  Fix a finite set \(\Lambda \),
which we will use as an \emph{alphabet}. We use \(\Lambda ^*\) to denote the
set of strings of elements of \(\Lambda \). A \emph{language} is any subset
\(L\subseteq \Lambda ^*\). We use \(\epsilon \) to denote the empty string.
Given \(\sigma \in\Lambda ^*\), \(\abs\sigma \) denotes the length of
\(\sigma \).
\begin{definition}
  The set of \emph{regular} languages over \(\Lambda \) is the smallest set of
  languages that contains the finite languages and is such that if
  \(L_1,L_2\subseteq \Lambda ^*\) are regular then so are \(L_1\cup L_2\),
  \(L_1L_2 = \set{\sigma \tau :\sigma \in L_1,\tau \in L_2}\), and \(L_1^* =
  \set{\sigma _1\cdots \sigma _n : \sigma _1,\ldots ,\sigma _n\in L_1}\).
\end{definition}
Regular languages can be characterized as those recognized by machines of a
certain form:
\begin{definition}
  A \emph{non-deterministic finite automaton} (or \emph{NFA}) is a \(5\)-tuple
  \(M = (\Lambda , Q,q_0,\Omega ,\delta )\), where:
  \begin{itemize}
    \item
      \(\Lambda \) is a finite alphabet.
    \item
      \(Q\) is a finite set of \emph{states}.
    \item
      \(q_0\in Q\) is the \emph{initial state}.
    \item
      \(\Omega \subseteq Q\) is the set of \emph{finish states}.
    \item
      \(\delta \colon Q\times \Lambda \to 2^Q\) is the \emph{transition
      function}: given \((q,a)\in Q\times \Lambda \) it outputs the set of
      states the machine could transition to if it is in state \(q\) and reads
      input \(a\).
  \end{itemize}
  We identify \(\delta \) with its natural extension \(Q\times \Lambda ^*\to
  2^Q\) given by \(\delta (q,\sigma a) = \bigcup \set{\delta (q',a) : q'\in
  \delta (q,\sigma )}\).  If \(\sigma \in\Lambda ^*\) we say \(M\)
  \emph{accepts} \(\sigma \) if \(\delta (q_0,\sigma )\cap \Omega \ne\emptyset
  \); the set of such \(\sigma \) is the language \emph{recognized} by \(M\).
  If \(\abs{\delta (q,a)} = 1\) for all \(q\in Q\) and \(a\in\Lambda \) we say
  \(M\) is a \emph{deterministic finite automaton} (or \emph{DFA}).  
\end{definition}
\begin{fact}[\cite[Lemma 2.2, Sections 3.2 and 3.3]{yu97}]
  If \(L\subseteq \Lambda ^*\) then \tfae{}:
  \begin{itemize}
    \item
      \(L\) is regular.
    \item
      \(L\) is recognized by some DFA.
    \item
      \(L\) is recognized by some NFA.
  \end{itemize}
\end{fact}
A note on exponential notation: we use \(\Lambda ^r\) to denote the
alphabet of \(r\)-tuples of elements of \(\Lambda \), and we use \(\Lambda
^{(r)}\) to denote the language \(\set{\sigma \in\Lambda ^* : \abs\sigma =r}\).

Throughout the paper, we fix an infinite abelian group \(\Gamma \) equipped
with some injective endomorphism \(F\colon \Gamma \to\Gamma \).  We let
\(\mathbb{Z}[F]\) denote the subring of \(\ndm(\Gamma )\) generated by \(F\),
and consider \(\Gamma \) as a \(\mathbb{Z}[F]\)-module. Our context is slightly
more general than that of \cite{bell19}: they require that \(\Gamma \) be a
finitely generated abelian group, which we do not. Most of the results of
\cite{bell19} go through in this context with no additional effort.

Following \cite{bell19}, given a string \(s_0\cdots s_n\) of elements in
\(\Gamma \) we let
\[[s_0\cdots s_n]_F = s_0 + F s_1 + \cdots + F^ns_n
\]
Note that when \((\Gamma ,F) = (\mathbb{Z},d)\) and \(s_i\in \set{0,\ldots
,d-1}\) this is just computing the number represented by \(s_0\cdots s_n\) base
\(d\). Given a set \(L\) of strings of elements of \(\Gamma \), we let
\[[L]_F = \set{[\sigma ]_F : \sigma \in L}.
\]
\begin{definition}
  \label{def:spanningsets}
  If \(\Sigma \) is a finite subset of \(\Gamma \) we say \(\Sigma \) is an
  \emph{\(F\)-spanning set} (for \(\Gamma \)) if it satisfies the following
  axioms:
  \begin{enumerate}[label=(\roman*)]
    \item
      For all \(a\in\Gamma \) there is \(\sigma \in\Sigma^* \) such that
      \(a = [\sigma ]_F\).
    \item
      \(0\in\Sigma \), and if \(a\in\Sigma \) then \(-a\in\Sigma \).
    \item
      If \(a_1,\ldots ,a_5\in\Sigma \) then \(a_1+\cdots +a_5\in \Sigma
      +F\Sigma \).
    \item
      If \(a_1,a_2,a_3\in\Sigma \) and \(a_1+a_2+a_3\in F\Gamma \) then
      \(a_1+a_2+a_3\in F\Sigma \).
  \end{enumerate}
\end{definition}
Conditions (ii-iv) are included largely for bookkeeping reasons; the main
condition is (i). It says that every element of \(\Gamma \) has an
``\(F\)-expansion with digits in \(\Sigma \)''. Note that the existence of an
\(F\)-spanning set will imply that \(\Gamma \) is finitely generated as a
\(\mathbb{Z}[F]\)-module, and that \(\Gamma /F\Gamma \) is finite.

It is pointed out in \cite[Lemmas 5.6 and 5.7]{bell19} that for \(r>0\) if
\(\Sigma \) is an \(F\)-spanning set then \([\Sigma ^{(r)}]_F\) is both an
\(F\)-spanning set and an \(F^r\)-spanning set.\footnote{What we denote
  \([\Sigma ^{(r)}]_F\) is denoted simply \(\Sigma ^{(r)}\) in \cite{bell19};
recall that for us \(\Sigma ^{(r)}\) denotes the set of strings over \(\Sigma
\) of length \(r\).}
\begin{definition}
  We say \(A\subseteq \Gamma \) is \emph{\(F\)-automatic} if there is an
  \(F^r\)-spanning set \(\Sigma \) for some \(r>0\) such that \(\set{\sigma
  \in\Sigma ^* : [\sigma ]_{F^r} \in A}\) is regular.
\end{definition}
The following is a useful strengthening of \cite[Proposition 6.3]{bell19}:
\begin{proposition}
  \label{prop:autoinspan}
  Suppose \(A\subseteq \Gamma \) is \(F\)-automatic. Then for any \(r>0\) and
  any \(F^r\)-spanning set \(\Sigma \), \(\set{\sigma \in\Sigma ^*:[\sigma
  ]_{F^r}\in A}\) is regular.
\end{proposition}
\begin{rpf}
  By \(F\)-automaticity there is \(r_0>0\) and an \(F^{r_0}\)-spanning set
  \(\Sigma _0\) such that \(\set{\sigma \in\Sigma _0^* : [\sigma ]_{F^{r_0}}\in
  A}\) is regular. Now \(\Theta := [\Sigma ^{(r_0)}]_{F^r}\) is an
  \(F^{rr_0}\)-spanning set, and by \cite[Proposition 6.3]{bell19} we have that
  \(\set{\sigma \in\Theta ^* : [\sigma ]_{F^{rr_0}}\in A}\) is regular. Suppose
  it is recognized by the automaton \(M=(\Theta ,Q,q_0,\Omega ,\delta )\). Now
  define a new automaton \(M' = (\Sigma ,Q\times \Sigma
  ^{(<r_0)},(q_0,\epsilon),\Omega ',\delta ')\) by
  \begin{ea}
    \delta '((q,\sigma ), a) &=& 
    \begin{cases}
      (\delta (q,[\sigma a]_{F^r}),\epsilon ) &\tif \abs{\sigma } = r_0-1 \\
      (q,\sigma a) &\tif \abs{\sigma }<r_0-1
    \end{cases} \\
    \Omega ' &=& \set{(q,\sigma ) : \delta (q,[\sigma ]_{F^r}) \in \Omega }
  \end{ea}
  (Note if \(\abs\sigma <r_0\) that \([\sigma ]_{F^r} = [\sigma
  0^{r_0-\abs\sigma }]_{F^r} \in \Theta \), so \(\delta (q,[\sigma ]_{F^r})\)
  is defined.) Given \(\sigma \in\Sigma ^*\) we can write \(\sigma =\sigma
  _1\cdots \sigma _{n+1}\) where \(\abs{\sigma _1}=\cdots =\abs{\sigma _n} =
  r_0\) and \(\abs{\sigma _{n+1}}< r_0\). Now, \(M'\) accepts \(\sigma \) if
  and only if \(M\) accepts \([\sigma _1]_{F^r}\cdots [\sigma _{n+1}]_{F^r}\)
  (where again \([\sigma _{n+1}]_{F^r} = [\sigma _{n+1}0^{r_0-\abs{\sigma
  _{n+1}}}]_{F^r}\in\Theta \)), which is in turn equivalent to \([\sigma
  ]_{F^r} = [[\sigma _1]_{F^r}\cdots [\sigma _{n+1}]_{F^r}]_{F^{rr_0}}\in A\).
  So \(M'\) recognizes \(\set{\sigma \in\Sigma ^*:[\sigma ]_{F^r}\in A}\), as
  desired.
\end{rpf}

\section{Existence of spanning sets}
We first study what it means for \(\Gamma \) to admit an \(F\)-spanning set. It
is clear in \cite{bell19} that the authors do not expect all finitely generated
abelian groups to admit spanning sets, but no example was given. Here is one:
\begin{example}
  \label{eg:luke}
  Consider \(\Gamma =\mathbb{Z}^2\).  Fix \(T\in M_2(\mathbb{Z})\) invertible
  and diagonalizable over \(\mathbb{C}\), and with an eigenvalue \(\mu \) with
  \(\abs\mu <1\); for example \(T = 
  \begin{pmatrix}
    1 & 1 \\
    1 & 0
  \end{pmatrix}\). Let \(F\colon \mathbb{Z}^2\to\mathbb{Z}^2\) be the
  associated linear map on \(\mathbb{Z}^2\).  Suppose for contradiction that
  there were an \(F\)-spanning set \(\Sigma \) for \(\mathbb{Z}^2 \). Let
  \(\set{v,w}\subseteq \mathbb{C}^2\) be an eigenbasis for \(T\), say with \(Tv
  = \mu v\) and \(Tw = \nu w\). Write each element of \(\Sigma \) (uniquely) in
  the form \(av+bw\) for \(a,b\in\mathbb{C}\), and let \(M\) be the largest
  \(\abs a\) obtained in this way.  Given \(x\in\mathbb{Z}^2\) we can write
  \[x = [s_0\cdots s_n]_F = \sum_{i=0}^n T^i(a_iv +
    b_iw) = \pars*{\sum_{i=0}^n a_i\mu ^i}v + \pars*{\sum_{i=0}^n b_i\nu ^i}w
  \]
  where \(s_i = a_i v+ b_iw\in \Sigma \). So by independence of \(\set{v,w}\)
  if we write \(x = av + bw\) then
  \[\abs a = \abs*{\sum_{i=0}^n a_i\mu ^i}
    \le \sum_{i=0}^n \abs{a_i} \abs\mu ^i
    \le M\sum_{i=0}^\infty \abs\mu ^i \le M(1-\abs\mu )^{-1}
  \]
  On the other hand, if we had \(a\ne 0\) then since \(\mathbb{Z}^2\) is closed
  under integer scaling there would be some \(x\in\mathbb{Z}^2\) such that the
  corresponding \(a\) did not satisfy \(\abs a\le M(1-\abs\mu )^{-1}\), a
  contradiction.  So \(a = 0\), and \(\mathbb{Z}^2\subseteq \mathbb{C}w\),
  contradicting the fact that the \(\mathbb{C}\)-linear span of
  \(\mathbb{Z}^2\) is all of \(\mathbb{C}^2\). So no \(F\)-spanning set exists
  for \(\mathbb{Z}^2\).
\end{example}
Bell and Moosa \cite{bell19} give a sufficient condition for \(\Gamma \) to
admit an \(F^r\)-spanning set for some \(r\) (under the assumption that
\(\Gamma /F\Gamma \) is finite) in terms of the existence of what they call a
\emph{height function} on \(\Gamma \). Our first goal is to show that this is
also necessary.
\begin{definition}
  A \emph{height function} for \((\Gamma,F) \) is a map \(h\colon \Gamma
  \to\mathbb{R}_{\ge0}\) satisfying the following:
  \begin{description}[font=\normalfont]
    \lit{(Symmetry and triangle inequality)}
      There are \(\alpha ,\kappa \in\mathbb{R}\) with \(D>1\) and \(\kappa >0\)
      such that \(h(-a)\le \alpha Dh(a)+\kappa \) and \(h(a+b)\le \alpha
      (h(a)+h(b))+\kappa \) for all \(a,b\in\Gamma \).
    \lit{(Northcott property)}
      For all \(N\in\mathbb{N}\) there are finitely many \(a\in\Gamma \) with
      \(h(a)\le N\).
    \lit{(Canonicity)}
      There is \(\beta \in\mathbb{R}\) with \(\beta >1\) such that
      \(h(Fa)>\beta h(a)\) for cofinitely many \(a\in\Gamma \).
  \end{description}
\end{definition}
We find it more useful to work with a variant of height that we will call
``length''.
\begin{definition}
  A \emph{length function} for \((\Gamma,F) \) is a map \(\lambda \colon
  \Gamma \to \mathbb{R}_{\ge0}\) satisfying the following:
  \begin{description}[font=\normalfont]
    \lit{(Symmetry)}
      \(\lambda (a) = \lambda (-a)\) for all \(a\in\Gamma \).
    \lit{(Ultrametric inequality)}
      There is \(D\in\mathbb{R}\) with \(D\ge 1\) such that \(\lambda (a+b)\le
      D\max(\lambda (a),\lambda (b))\) for all \(a,b\in\Gamma \).
    \lit{(Northcott property)}
      For all \(N\in\mathbb{N}\) there are finitely many \(a\in\Gamma \) with
      \(\lambda (a)\le N\).
    \lit{(Canonicity)}
      There is a finite exceptional set \(A\) and \(C,E\in\mathbb{R}\) with
      \(C>1\) and \(E\ge1\) such that
      \begin{itemize}
        \item
          \(\lambda (Fa)\le C\lambda (a)\) for all \(a\in\Gamma \), and
        \item
          \(\lambda (F^na)\ge \frac{C^n}E\lambda (a)\) for all \(a\in\Gamma
          \setminus A\) and \(n\in\mathbb{N}\).
      \end{itemize}
  \end{description}
\end{definition}
Note that if \(\lambda \) satisfies all the axioms of being a length function
besides symmetry, then \(\lambda '(a) := \max(\lambda (a) , \lambda (-a))\)
will be a length function.

Length and height functions are closely related. Suppose \(\lambda \) is a
length function for \((\Gamma ,F)\) with exceptional set \(A\) and associated
constants \(C,D,E\). Pick \(r\) so that \(C^r>E\). It is then straightforward
to check that \(\lambda \) is a height function for \((\Gamma ,F^r)\) with
exceptional set \(A\) and associated constants \(\alpha =D\), \(\beta =
C^rE^{-1}\), and \(\kappa =0\).  It is harder to derive a length function from
a height function, as the axioms of height functions place no upper bound on
the height of \(Fa\) in terms of the height of \(a\).  However, the height
functions that arise in the isotrivial Mordell-Lang context of \cite{bell19}
turn out to also be length functions.
\begin{remark}
  \label{rem:lengthlargec}
  It will sometimes be convenient to assume that \(C\) is large compared to
  some function of \(D\) and \(E\) (and possibly other constants). If we are
  willing to pass from \(F\) to a power thereof, this can always be assumed: if
  \(\lambda \) is a length function for \((\Gamma ,F)\) with associated
  constants \(C,D,E\), we can take \(r\) such that \(C^r\) satisfied the
  desired inequalities. Then \(\lambda \) is a length function for \((\Gamma
  ,F^r)\) with associated constants \(C^r,D,E\).
\end{remark}
\begin{lemma}
  \label{lemma:canonicityinforb}
  Suppose \(\lambda \) is a length function for \((\Gamma ,F)\) with
  exceptional set \(A\) and associated constants \(C,D,E\).  By increasing
  \(E\) we can assume that every element of the exceptional set \(A\) has
  finite \(F\)-orbit.
\end{lemma}
\begin{rpf}
  Suppose the \(F\)-orbit of \(a\in A\) is infinite; so there is \(i_a\) such
  that \(F^{i_a}a\notin A\). Note that all \(\lambda (F^ia) \ne 0\): otherwise
  the ultrametric inequality would imply \(\lambda (kF^ia) = 0\) for all \(k\),
  contradicting the Northcott property.  Now, for all \(n\ge i_a\) we have
  \[\lambda (F^na) = \lambda (F^{n-i_a}F^{i_a}a) \ge \frac{C^{n-i_a}}E \lambda
    (F^{i_a}a) \ge C^n\lambda (a) \pars*{\underbrace{C^{i_a}E\frac{\lambda
    (a)}{\lambda (F^{i_a}a)}}_*}^{-1}
  \]
  and for \(n< i_a\) we have
  \[\lambda (F^na) = C^n\lambda (a) \pars*{\underbrace{C^n\frac{\lambda
    (a)}{\lambda (F^na)}}_*}^{-1}
  \]
  Hence taking \(E'\) to be the maximum of \(E\) and all the quantities marked
  \(*\), we get that \(\lambda (F^na)\ge \frac{C^n}{E'}\lambda (a)\) so that
  \(a\) is no longer exceptional with respect to this new \(E'\). Iterating
  this procedure for each \(a\in A\) with infinite \(F\)-orbit, we eventually
  produce a new \(\widetilde{E}\ge E\) such that only the elements of \(A\) of
  finite \(F\)-orbit remain exceptional.
\end{rpf}
\begin{remark}
  \label{rem:inforderinforbit}
  If there is a length function for \((\Gamma ,F)\) and \(a\in\Gamma \) has
  finite \(F\)-orbit, then \(a\) must have finite order.  Indeed, the
  \(F\)-orbit of \(ka\) is finite for all \(k\in\mathbb{N}\); but canonicity
  then implies that all \(ka\) lie in the finite exceptional set.
\end{remark}
Combining \cref{lemma:canonicityinforb,rem:inforderinforbit} we may assume
canonicity applies to all elements of infinite order.

Here is the motivating example of a length function.
\begin{definition}
  \label{def:spanlenfn}
  If \(\Sigma \) is an \(F\)-spanning set for \(\Gamma \), we define \(\lambda
  _\Sigma \colon \Gamma \to\mathbb{R}_{\ge0}\) by \(\lambda _\Sigma (a) =
  2^\ell \), where \(\ell \) is the length of the shortest \(\sigma \in\Sigma
  ^*\) such that \([\sigma ]_F = a\).
\end{definition}
\begin{proposition}
  \label{prop:lenabstractlen}
  If \(\Sigma \) is an \(F\)-spanning set for \(\Gamma \) then \(\lambda
  _\Sigma \) is a length function for \((\Gamma ,F)\) with associated constants
  \(C=D=E=2\) and exceptional set \(\Sigma \).
\end{proposition}
\begin{rpf}
  We verify the axioms.  Symmetry of \(\lambda _\Sigma \) is simply by symmetry
  of \(\Sigma \). Lemma 5.3 of \cite{bell19} implies that the ultrametric
  inequality holds of \(\lambda _\Sigma \) with \(D=2\). Since \(\Sigma \) is
  finite there are finitely many \(\sigma \in\Sigma ^*\) of length \(\le N\),
  and hence \(\lambda _\Sigma \) satisfies the Northcott property. Finally, we
  prove canonicity for \(\lambda _\Sigma \) with \(C=E=2\) and exceptional set
  \(\Sigma \). For the upper bound, note that if \(a = [s_0\cdots s_\ell ]_F\)
  with \(s_0,\ldots ,s_\ell \in\Sigma \) then \(Fa = [0s_0\cdots s_\ell ]_F\);
  so \(\lambda _\Sigma (Fa)\le 2\lambda _\Sigma (a)\). It remains to show the
  lower bound.

  Suppose \(a\notin \Sigma \). Write \(\lambda _\Sigma (a) = 2^\ell \) for some
  \(\ell >1\); say \(a = [s_0\cdots s_{\ell-1}]_F\). Suppose for contradiction
  that \(\lambda _\Sigma (F^ma) < 2^{m+\ell -1} \); so we can write \(F^ma =
  [t_0\cdots t_{m+\ell -3}]_F\). Then
  \begin{equation}
    \label{eqn:canonicitysource}
    F^ms_0+\cdots + F^{m+\ell -1}s_{\ell -1} = F^ma = t_0 + \cdots +
    F^{m+\ell -3}t_{m+\ell -3}
  \end{equation}
  Then \(t_0\in F\Gamma \), so by axiom (iv) we get that \(t_0 = Ft_0'\) for
  some \(t_0'\in \Sigma \).  Inductively suppose for some \(i<m-1\) we can
  write \(t_0 + \cdots + F^it_i = F^{i+1}t_i'\) for some \(t_i'\in\Sigma \).
  Then \(F^{i+1}(t_i'+t_{i+1}) = t_0 + \cdots + F^{i+1}t_{i+1}\) which can be
  seen to be in \( F^{i+2}\Gamma \) using \cref{eqn:canonicitysource} and the
  fact that \(m\ge i+2\). Hence \(t_i'+t_{i+1}\in F\Gamma \) by injectivity
  of \(F\), and again by axiom (iv) we can write \(t_i'+t_{i+1} = Ft_{i+1}'\)
  for some \(t_{i+1}\in\Sigma \). So \(t_0 + \cdots + F^{i+1}t_{i+1} =
  F^{i+2}t_{i+1}'\).  It follows that \(t_0+\cdots +F^{m-1}t_{m-1} =
  F^mt_{m-1}'\) for some \(t_{m-1}'\in\Sigma \). (Note \(t_{m-1}\) is defined
  since \(\ell \ge2\).) So \(F^ma = F^mt_{m-1}' + F^mt_m + \cdots + F^{m+\ell
  -3}t_{m+\ell -3}\), and thus
  \[a = t_{m-1}' + (t_m + \cdots + F^{\ell -3}t_{m+\ell -3})
  \]
  Hence by \cite[Lemma 5.3]{bell19} we get that \(a\) can be represented by
  a string of length \(\le \ell -1\), contradicting our assumption that
  \(\lambda _\Sigma (a) = 2^\ell \).
\end{rpf}
Putting all this together, we deduce the following characterization for the
existence of spanning sets, which in particular provides a converse to
\cite[Proposition 5.8]{bell19}.
\begin{theorem}
  \label{thm:spanlenheight}
  Suppose \(\Gamma /F\Gamma \) is finite.
  \Tfae{}:
  \begin{enumerate}
    \item
      \(\Gamma \) admits an \(F^r\)-spanning set for some \(r>0\).
    \item
      There is a length function for some \((\Gamma ,F^r)\).
    \item
      There is a height function for some \((\Gamma ,F^r)\).
  \end{enumerate}
\end{theorem}
\begin{rpf}
  (1)\(\implies \)(2) is the previous proposition, and we noted (2)\(\implies
  \)(3) after the definition of length functions.  For (3)\(\implies \)(1), we
  appeal to \cite[Proposition 5.8]{bell19}. Formally they require that \(\Gamma
  \) be finitely generated as a group, not merely as a
  \(\mathbb{Z}[F]\)-module, but this is only used to deduce that \(\Gamma
  /F^r\Gamma \) is finite for all \(r\); we show that we can assume as much.

  Let \(S\subseteq \Gamma \) contain a representative of each coset of
  \(F\Gamma \); we claim that \([S^{(r)}]_F = \set{s_0 + \cdots +
  F^{r-1}s_{r-1} : s_i\in S}\) contains a representative of each coset of
  \(F^r\Gamma \). Indeed, given \(a\in\Gamma \) we can find \(s_0\in S\) such
  that \(a\equiv s_0\pmod{F\Gamma }\). Then inductively we can find
  \(s_1,\ldots ,s_{r-1}\in S\) such that \(F^{-1}(a-s_0)\equiv s_1+\cdots
  +F^{r-2}s_{r-1}\pmod{F^{r-1}\Gamma }\), at which point it follows that
  \(a\equiv s_0+\cdots +F^{r-1}s_{r-1}\pmod{F^r\Gamma }\). 

  Since \(\Gamma /F\Gamma \) is finite, we can take \(S\) to be finite. So
  \(\abs{\Gamma /F^r\Gamma }\le \abs{[S^{(r)}]_F}\le \abs{S}^r<\infty \), and
  \(\Gamma /F^r\Gamma \) is finite, as desired.
\end{rpf}
\begin{corollary}
  \label{cor:spanfinvsub}
  Suppose there is an \(F^r\)-spanning set for \(\Gamma \) for some \(r>0\) and
  \(H\le \Gamma \) is \(F\)-invariant. Then there is an \(F^s\)-spanning set
  for \(H\) for some \(s>0\). Furthermore if \(A\subseteq H\) then \(A\) is
  \(F\)-automatic in \(\Gamma \) if and only \(A\) is \(F\)-automatic in \(H\).
\end{corollary}
\begin{rpf}
  By \cref{thm:spanlenheight} there is a length function \(\lambda \) for some
  \((\Gamma ,F^r)\). One can check that \(\lambda \restriction H\) is a length
  function for \((H,(F\restriction H)^r)\). So again by
  \cref{thm:spanlenheight} there is an \(F^s\)-spanning set for \(H\) for some
  \(s>0\).

  For the ``furthermore'', note by the above and \cite[Lemmas 5.6 and
  5.7]{bell19} there is \(s\) such that there is an \(F^s\)-spanning set
  \(\Sigma \) for \(H \) and that is contained in an \(F^s\)-spanning set
  \(\Sigma '\) for \(\Gamma \). The right-to-left direction then follows from
  \cite[Proposition 6.8 (b)]{bell19}. For the left-to-right, note that
  \(\set{\sigma \in \Sigma ^* : [\sigma ]_{F^s}\in A} = \set{\sigma \in(\Sigma
  ')^* : [\sigma ]_{F^s}\in A}\cap\Sigma ^*\) is the intersection of two
  regular languages, and is thus regular (see \cite[Theorem 2.1]{yu97}); so
  \(A\) is \(F\)-automatic in \(H\).
\end{rpf}
In fact when \(\Gamma \) is a finitely generated abelian group, we can use the
above to deduce a concrete and verifiable characterization of the existence of
spanning sets.
\begin{theorem}
  \label{thm:spaneigen}
  Suppose \(\Gamma \) is a finitely generated abelian group. Then \(\Gamma \)
  admits an \(F^r\)-spanning set for some \(r>0\) if and only if the
  eigenvalues of \(F\otimes _\mathbb{Z}\id_\mathbb{C}\) (viewed as a linear map
  on the \(\mathbb{C}\)-vector space \(\Gamma \otimes _\mathbb{Z}\mathbb{C}\))
  all have modulus \(>1\).
\end{theorem}
\begin{rpf}
  We first show that we can reduce to the torsion-free case. By the fundamental
  theorem of finitely generated abelian groups we may assume \(\Gamma
  =\mathbb{Z}^m \times H\) where \(H\) is a finite group. Since \(F\) is
  injective we get that \(F^{-1}H = H\); so \(F\) induces an injective
  endomorphism \(F'\colon \mathbb{Z}^m\to\mathbb{Z}^m\) on the quotient
  \(\Gamma /H\cong\mathbb{Z}^m\).  Then \(\Gamma \) has an \(F\)-spanning set
  if and only if \((\mathbb{Z}^m,F')\) does: it is not hard to verify that if
  \(\Sigma \) is an \(F\)-spanning set for \(\Gamma \) and \(\pi \colon \Gamma
  \to \mathbb{Z}^m\) is the projection then \(\pi (\Sigma)\) is an
  \(F'\)-spanning set for \(\mathbb{Z}^m\), and conversely if \(\Sigma \) is an
  \(F'\)-spanning set for \(\mathbb{Z}^m\) then \(\pi ^{-1}(\Sigma) \) is an
  \(F\)-spanning set for \(\Gamma \).  Note further that under the
  identification \(\Gamma \otimes _\mathbb{Z}\mathbb{C}=\mathbb{Z}^m\otimes
  _\mathbb{Z}\mathbb{C}\) we have that \(F\otimes _\mathbb{Z}\id_\mathbb{C}
  \sim F'\otimes _\mathbb{Z}\id_\mathbb{C}\), and in particular they have the
  same eigenvalues.

  We may therefore assume \(\Gamma =\mathbb{Z}^m\) for some \(m\) and
  \(F\in M_m(\mathbb{Z})\).
  \begin{description}
    \prl By replacing \(F\) with a power, we may assume \(\abs\mu  >2\) for all
      \(\mu \in\sigma (F)\). We show there exists a height function for
      \((\Gamma ,F)\).  Pick a basis \(\set{e_1,\ldots ,e_m }\) for
      \(\mathbb{C}^m\) that puts \(F\) in Jordan canonical form. Let \(h\colon
      \mathbb{C}^m\to\mathbb{R}_{\ge0}\) be the infinity norm associated to
      \(\set{e_1,\ldots ,e_m}\): so \(h(a_1e_1+\cdots +a_me_m) = \max_i
      \abs{a_i}\). We show that \(h\restriction \mathbb{Z}^m\) is a height
      function for \((\Gamma ,F)\). Triangle inequality and symmetry are clear;
      Northcott follows from the fact that all norms on a finite-dimensional
      space are equivalent (and in particular that our infinity norm is
      equivalent to the usual infinity norm on \(\mathbb{C}^m\)). It remains to
      check canonicity.

      Fix some \(\beta >1\) such that \(\beta +1<\abs\mu \) for all \(\mu
      \in\sigma (F)\). Suppose \(v = v_1e_1+\cdots
      +v_me_m\in\mathbb{C}^m\setminus \set0\); write \(Fv = w_1e_1+\cdots
      +w_me_m\). Fix \(i\) such that \(\abs{v_i} = h(v)\).  Depending on the
      structure of the Jordan blocks, we get that \(w_i\) is either \(\mu v_i +
      v_{i+1}\) or \(\mu v_i\) (for \(\mu \in\sigma (F)\) corresponding to
      \(v_i\)). In the first case, reverse triangle inequality yields
      \[h(Fv)\ge \abs{w_i}\ge \abs{\mu }\abs{v_i} - \abs{v_{i+1}}
        \ge (\abs\mu -1)\abs{v_i} > \beta \abs{v_i} = \beta h(v)
      \]
      and in the second case we get
      \[h(Fv)\ge \abs{w_i} = \abs\mu \abs{v_i} > \beta \abs{v_i} = \beta h(v)
      \]
      So \(h\) satisfies canonicity, as desired.
    \plr
      Suppose there is \(\mu \in\sigma (F)\) with \(\abs\mu \le1\). There are
      two cases:
      \begin{caselist}
        \case
          Suppose there is \(\mu \in\sigma (F)\) with \(\abs\mu <1\); the
          argument in this case is similar to \cref{eg:luke}. Since the same is
          true for all powers of \(F\) it suffices to show that there is no
          \(F\)-spanning set. Pick a basis \(\set{e_1,\ldots ,e_m}\) for
          \(\mathbb{C}^m\) that puts \(F\) into Jordan canonical form; for
          \(x\in\mathbb{C}^m\) we write
          \[x = \sum_{i=1}^m f_i(x)e_i
          \]
          so each \(f_i\colon \mathbb{C}^m\to\mathbb{C}\) is linear. Pick some
          \(e_k\) corresponding to the bottom-right of some Jordan block for
          \(\mu \); so \(f_k(Fx) = \mu f_k(x)\) for \(x\in\mathbb{C}^m\).

          Suppose for contradiction we had an \(F\)-spanning set \(\Sigma \).
          Let \(M = \max\set{\abs{f_k(a)}:a\in\Sigma }\). Then if
          \[y = \sum_{i=0}^\ell F^ia_i
          \]
          for \(a_i\in\Sigma \) then
          \[\abs{f_k(y)} = \abs*{\sum_{i=0}^\ell \mu ^if_k(a_i)}
            \le \sum_{i=0}^\ell \abs{\mu ^if_k(a_i)}
            = \sum_{i=0}^\ell \abs\mu ^i\abs{f_k(a_i)}
            \le M\sum_{i<\omega }\abs\mu ^i
            =M\frac1{1-\abs\mu }
          \]
          So every \(y\in\mathbb{Z}^m\) satisfies \(\abs{f_k(y)}\le
          M\frac1{1-\abs\mu }\). So since the integers are closed under
          doubling and \(f_k\) is linear, we get that \(f_k(y) = 0\) for each
          \(y\in\mathbb{Z}^m\). So \(\mathbb{C}^m =
          \lspan_\mathbb{C}\mathbb{Z}^m\) is spanned by
          \(\set{e_1,\ldots,e_m}\setminus \set{e_k}\), a contradiction.
        \case
          Suppose \(\abs\mu \ge1\) for all \(\mu \in\sigma (F)\); pick some
          \(\mu \in\sigma (F)\) with \(\abs\mu =1\). We show that there is
          non-zero \(a\in\mathbb{Z}^m\) such that \(F^ia=F^ja\) for some \(i\ne
          j\), and deduce that no power of \(F\) admits a length function.

          Let \(p_\mu (x)\) be the minimal polynomial for \(\mu \) over
          \(\mathbb{Q}\). Note that if \(\nu \in\mathbb{C}\) is a root of
          \(p_\mu \) then so is \(\nu ^{-1}\).  Indeed, this is true of \(\mu
          \) since \(\mu ^{-1}=\overline{\mu }\) and \(p_\mu
          \in\mathbb{Q}[x]\), and hence is true of all
          \(\aut(\overline{\mathbb{Q}}/\mathbb{Q})\)-conjugates of \(\mu \).
          But \(p_\mu \) divides the characteristic polynomial of \(F\), and
          hence \(p_\mu \) has no roots of modulus \(<1\) by hypothesis. It
          follows that \(p_\mu \) has no roots of modulus \(>1\) either. That
          is, the roots of \(p_\mu \) lie on the unit circle.

          Now, \(p_\mu (F)\) isn't invertible over \(\mathbb{C}\), and hence
          isn't invertible over \(\mathbb{Q}\); so \(V:=\ker_\mathbb{Q}(p_\mu
          (F))\) is a non-trivial \(F\)-invariant subspace of \(\mathbb{Q}^m\).
          Furthermore \(p_\mu (F\restriction V) = 0\), so the minimal
          polynomial of \(F\restriction V\) is separable (since \(p_\mu \) is);
          so \(F\restriction V\) is diagonalizable (over \(\mathbb{C}\)) with
          eigenvalues that are roots of \(p_\mu \), and hence lie on the unit
          circle. In particular in the inner product induced by the eigenbasis
          for \(F\restriction V\) we get that \(F\restriction V\) is unitary.

          Since \(V\) is a non-trivial subspace of \(\mathbb{Q}^m\) it contains
          a non-zero \(a\in\mathbb{Z}^m\). So \(Fa\in V\cap\mathbb{Z}^m\) as
          well and \(\nrm{Fa} = \nrm{a}\) (where the norm is induced by the
          aforementioned inner product); inductively we get \(\nrm{F^na} = \nrm
          a\) and \(F^na\in\mathbb{Z}^m\) for all \(n\).  But \(\nrm{\cdot }\)
          is equivalent to the Euclidean norm (since both are norms on a
          finite-dimensional space); so there are finitely many \(b\in
          \mathbb{Z}^m\) with \(\nrm{b} = \nrm{a}\).  So \(F^i a = F^ja\) for
          some \(i\ne j\). By injectivity of \(F\) we get \(F^{i-j}a = a\), and
          hence that \((F^n)^{i-j} a = a\) for all \(n\). It follows that the
          \(F^n\)-orbit of \(a\) is finite for all \(n>0\). Since \(a\) has
          infinite order, \cref{rem:inforderinforbit} yields that no power of
          \(F\) admits a length function. Hence by \cref{thm:spanlenheight}
          \(\Gamma \) does not admit an \(F^r\)-spanning set for any \(r>0\).
          \qedhere
      \end{caselist}
  \end{description}
\end{rpf}
\section{A characterization of \(F\)-automaticity}
A useful characterization of classical automaticity is in terms of ``finiteness
of kernels''. A version of this for \(F\)-automaticity is given in \cite[Lemma
6.2]{bell19}. Using length functions we are able to improve this result. First,
here is what kernels mean in our setting, as introduced in Definition 6.1 of 
\cite{bell19}. As before, we fix an abelian group \(\Gamma \) with an injective
endomorphism \(F\).
\begin{definition}
  Suppose \(A\subseteq \Gamma \). Given \(s_0,\ldots ,s_{n-1}\in\Gamma \) we
  set
  \[A_{s_0\cdots s_{n-1}}(x) =
    \set{x\in\Gamma  : s_0 + \cdots  + F^{n-1}s_{n-1}+F^n x\in A}
  \]
  (so \(A_\epsilon = A\)). Given \(S\subseteq \Gamma \), by the
  \emph{\((S,F)\)-kernel of \(A\)} we mean the set \(\set{A_\sigma :\sigma \in
  S^*} \).
\end{definition}
\begin{theorem}
  \label{thm:cosetdefauto}
  Suppose \(\Gamma \) admits an \(F^r\)-spanning set for some \(r>0\) (and so
  in particular that \(\Gamma /F\Gamma \) is finite). Fix finite \(S\subseteq
  \Gamma \) containing a representative of each coset of \(F\Gamma \) in
  \(\Gamma \). Then \(A\subseteq \Gamma \) is \(F\)-automatic if and only if
  the \((S,F)\)-kernel of \(A\) is finite.
\end{theorem}
The big improvement here is that \(S\) need not be a spanning set---being a
complete set of representatives for \(\Gamma /F\Gamma \) is a much weaker
condition.  A nice feature of this characterization is that, apart from the
hypothesis on existence of a spanning set, \(F\)-automaticity of a set \(A\)
can be checked without references to higher powers of \(F\).
\begin{lemma}
  Let \(\Gamma \), \(F\), and \(S\) be as in \cref{thm:cosetdefauto}. Suppose
  \(n>0\). The \((S,F)\)-kernel of \(A\) is finite if and only if the
  \((S^{(n)},F^n)\)-kernel of \(A \) is finite.
\end{lemma}
\begin{rpf}
  The \((S^{(n)},F^n)\)-kernel is clearly contained in the \((S,F)\)-kernel.
  Suppose the \((S,F)\)-kernel is infinite. 

  Note that if \(\sigma =s_0\cdots s_{k-1},\tau =t_0\cdots t_{\ell -1}\in S^*\)
  and \(a\in\Gamma \) then
  \begin{ea}
    a\in (A_\sigma )_\tau &\iff& t_0+\cdots +F^{\ell -1}t_{\ell -1}+F^\ell a\in
    A_\sigma \\
    &\iff& s_0 + \cdots + F^{k-1}s_{k-1} + F^k(t_0+\cdots +F^{\ell -1}t_{\ell
    -1}+F^\ell a) \in A \\
    &\iff& s_0 + \cdots + F^{k-1}s_{k-1} + F^kt_0 + \cdots + F^{k+\ell
    -1}t_{\ell -1} + F^{k+\ell }a \in A \\
    &\iff& a\in A_{\sigma \tau}
  \end{ea}
  So \((A_\sigma )_\tau  = A_{\sigma \tau }\); in particular, if \(A _{\sigma
  \tau }\ne A _{\sigma '\tau }\) then \(A _\sigma \ne A _{\sigma'}\).  

  Since the set of all \(A_\sigma \) is infinite, there is \(0\le i<n\) such
  that the set of \(A_\sigma \) with \(\abs\sigma \equiv i\pmod n\) is
  infinite. Now, as \(S\) is finite, this means there is \(\tau \in S^*\) of
  length \(i\) such that \(\set{A _{\rho \tau } : \rho \in S^*,\abs\rho \in
  n\mathbb{Z}} \) is infinite.  It follows by the above observation that
  \(\set{A _\rho  :\rho \in S^*,\abs\rho \in n\mathbb{Z}} \) is infinite, and
  so the \((S^{(n)}, F^n)\)-kernel is infinite.
\end{rpf}
\begin{lemma}
  Suppose \(\Gamma \) admits a length function \(\lambda \) for \((\Gamma ,F)\)
  with associated constants \(C,D,E\) such that \(C\ge D\).  Suppose
  \(S,T\subseteq \Gamma \) are finite sets both containing a representative of
  each coset of \(F\Gamma \).  Then \(A\) has finite \((S,F)\)-kernel if and
  only if it has finite \((T,F)\)-kernel.
\end{lemma}
\begin{rpf}
  Let \(N = \max\set{\lambda (s-t) : s\in S,t\in T}\). Suppose the
  \((T,F)\)-kernel of \(A\) is finite.
  
  Take an element \(A _\sigma  \) of the \((S,F)\)-kernel of \(A \), say with
  \(\sigma = s_0\cdots s_{n-1}\in S^*\). By the proof of
  \cref{thm:spanlenheight} there is \(\tau = t_0\ldots t_{n-1}\in T^*\) such
  that
  \[s_0 + \cdots + F^{n-1}s_{n-1}\equiv t_0+\cdots
    +F^{n-1}t_{n-1}\pmod{F^n(\Gamma )}
  \]
  So \(A_\tau \) lies in the \((T,F)\)-kernel of \(A \). Let
  \[\delta = F^{-n}((s_0-t_0) + \cdots + F^{n-1}(s_{n-1}-t_{n-1}))
  \]
  so for \(a\in\Gamma \) we have
  \begin{ea}
    a+\delta \in A_\tau  
    &\iff& t_0 + \cdots +F^{n-1}t_{n-1} + F^n (a+\delta )\in A\\
    &\iff& t_0 + \cdots +F^{n-1}t_{n-1} + F^n a + (s_0-t_0)+\cdots
    +F^{n-1}(s_{n-1}-t_{n-1}) \in A\\
    &\iff& s_0 + \cdots F^{n-1}s_{n-1} + F^na \in A\\
    &\iff& a\in A_\sigma 
  \end{ea}
  So \(A_\sigma = A_\tau -\delta \). Hence to show that the \((S,F)\)-kernel of
  \(A\) is finite it suffices to show that \(\delta \) can take on one of only
  finitely many values (as \(\sigma ,\tau \) vary). We show that \(\lambda
  (\delta )\) is bounded in terms of \(S,T\), which suffices by the Northcott
  property.

  We first bound \(\lambda (F^n\delta ) = \lambda ((s_0-t_0) + \cdots +
  F^{n-1}(s_{n-1}-t_{n-1})) \). Since \(\lambda (s_0-t_0)\le N\) and \(\lambda
  (F(s_1-t_1))\le CN\), the ultrametric inequality yields
  \[\lambda ((s_0-t_0) + F(s_1-t_1))\le DCN
  \]
  Since \(\lambda (F^2(s_2-t_2))\le C^2N\), another application of ultrametric
  inequality yields
  \[\lambda ((s_0-t_0) + F(s_1-t_1) + F^2(s_2-t_2))\le DC^2N
  \]
  (since by hypothesis \(C\ge D\)). Continuing inductively, we find that
  \(\lambda (F^n\delta )\le DC^{n-1}N\). Hence by canonicity we get that
  \(\lambda (\delta )\le EC^{-1}DN\).
\end{rpf}
\begin{lrpf}{thm:cosetdefauto}
  Suppose \(\Gamma \) admits an \(F^n\)-spanning set \(\Sigma \) for some
  \(n>0\). By \cref{prop:lenabstractlen} \(\lambda _\Sigma \)
  is a length function for \((\Gamma ,F^n)\) with associated constants
  \(C,D,E\) such that \(C\ge D\). Suppose \(A\subseteq \Gamma \). Then by the
  previous two lemmas, the \((S,F)\)-kernel of \(A\) is finite if and only if
  the \((S^{(n)},F^n)\)-kernel of \(A\) is, which occurs if and only if the
  \((\Sigma ,F^n)\)-kernel of \(A\) is. (Note that an \(F^n\)-spanning set must
  contain a representative of every coset of \(F\Gamma \).) But by \cite[Lemma
  6.2]{bell19} this is equivalent to \(A\) being \(F\)-automatic, as desired.
\end{lrpf}
As an illustration of the usefulness of this characterization we give quick
proof of \cite[Lemma 6.7(a)]{bell19}: that if a spanning set exists then the
equivalence relation of ``representing the same element'' on strings is
characterized by a finite automaton.
\begin{corollary}
  \label{cor:diagauto}
  Suppose \(\Sigma \) is an \(F\)-spanning set. Then \(\set{(\sigma ,\tau )\in
  (\Sigma ^2)^* : [\sigma ]_F = [\tau ]_F} \) is a regular language over
  \((\Sigma ^2)^*\).
\end{corollary}
\begin{rpf}
  By \cref{prop:autoinspan} this is equivalent to \(\Delta = \set{(a,a) :
  a\in\Gamma }\subseteq \Gamma ^2\) being \(F\)-automatic.  Fix a set \(S\)
  containing exactly one representative of each coset of \(F\Gamma \) in
  \(\Gamma \); we will show that the \((S^2,F)\)-kernel of \(\Delta \) is
  finite.  Note that if \(s_0\cdots s_{n-1},t_0\cdots t_{n-1}\in S^{(n)}\) are
  unequal, say with \(i\) minimal such that \(s_i\ne t_i\), then since
  \(s_i\centernot\equiv t_i\pmod{F\Gamma }\) and \(F\) is injective we get that
  \(s_0+\cdots +F^{n-1}s_{n-1}\centernot\equiv t_0+\cdots +
  F^{n-1}t_{n-1}\pmod{F^{i+1}\Gamma }\). So \(s_0+\cdots
  +F^{n-1}s_{n-1}\centernot\equiv t_0+\cdots +F^{n-1}t_{n-1}\pmod{F^n\Gamma
  }\), and thus \(s_0 + \cdots + F^{n-1}s_{n-1} + F^na \ne t_0 + \cdots +
  F^{n-1}t_{n-1} + F^nb\) for any \(a,b\in\Gamma \).  So \(\Delta _{(\sigma
  ,\tau) }\) is empty if \((\sigma, \tau) \in (S^2)^*\) and \(\sigma \ne\tau
  \). Furthermore if \(\sigma =\tau \) then by injectivity of \(F\) we get
  \(\Delta _{(\sigma ,\tau )} = \Delta \). So the \((S^2,F)\)-kernel of
  \(\Delta \) contains two elements.
\end{rpf}

\section{A characterization of \(F\)-sparsity}
\label{sec:sparsity}
Recall that a language \(L\subseteq \Sigma ^*\) is \emph{sparse} if it is
regular and \(\abs{\set{\sigma \in L : \abs\sigma \le x}}\) grows polynomially
in \(x\). See the beginning of \cite[Section 7]{bell19} for a brief overview
of sparsity.

Fix \((\Gamma ,F)\) an abelian group equipped with an injective endomorphism.
In \cite{bell19} Bell and Moosa adapt the notion of sparsity to the setting of
\(F\)-automatic sets. 
\begin{definition}
  \label{def:fsparse}
  A subset \(A\subseteq \Gamma \) is \emph{\(F\)-sparse} if there is an
  \(F^r\)-spanning set \(\Sigma \) for some \(r>0\) and a sparse \(L\subseteq
  \Sigma ^*\) such that \(A = [L]_{F^r}\).
\end{definition}
While the definition of \(F\)-sparsity is sufficient for the purposes of
\cite{bell19}, it is cumbersome in some contexts. In particular, in order to
show a set is not sparse one needs to check every possible set of
representatives in every possible spanning set for every possible power of
\(F\). For example, it is not immediate from the definition that \(\Gamma \)
itself isn't \(F\)-sparse.

In this section we will give a characterization of \(F\)-sparsity in terms of
length functions, and deduce from that some natural properties of
\(F\)-sparsity.  We will use the following characterization of sparsity:
\begin{fact}[\cite[Proposition 7.1]{bell19}]
  \label{fact:sparsechar}
  \(L\subseteq \Sigma ^*\) is sparse if and only if it is a finite union of
  sets of the form \(v_0w_1^*v_1\cdots v_{n-1}w_n^*v_n\) where \(v_0,\ldots
  ,v_n,w_1,\ldots ,w_n\in\Sigma ^*\).
\end{fact}
We call languages of the form \(v_0w_1^*v_1\cdots v_{n-1}w_n^*v_n\)
\emph{simple sparse}. By \cref{fact:sparsechar} an \(F\)-sparse set is a
finite union of sets of the form \([L]_{F^r}\) where \(L\) is a simple sparse
language. The following further simplification in the structure of \(F\)-sparse
sets will be useful both for proving basic closure properties below, and for
studying stability in the next section.
\begin{proposition}
  \label{prop:sparsenoconstants}
  Suppose \(A\subseteq \Gamma \) is \(F\)-sparse. There is \(s_0\in\mathbb{N}\)
  such that for all \(s\in s_0\mathbb{N}\) we can write \(A\) as a finite union
  of translates of sets of the form \([a _1^*\cdots a _n^*]_{F^s}\) where
  \(a_1,\ldots ,a_n\in \Gamma \).
\end{proposition}
\begin{rpf}
  By \cref{fact:sparsechar} we can write \(A\) as a finite union of sets of the
  form \([v_0w_1^*v_1\cdots w_n^*v_n]_{F^r}\) for some strings
  \(v_i,w_i\in\Gamma ^*\).  By replacing \(F\) with \(F^r\) we may assume
  \(r=1\). Let \(s_0\) be the least common multiple of all the \(\abs{w_i}\).

  In the case where \(\Gamma =\mathbb{Z}^m\) and \(F\) is multiplication by
  some \(d>0\), an intermediate result in the proof of \cite[Lemma
  3.5]{hawthorne20} is that we can write \(A\) as a finite union of translates
  of sets of the form \([\tau _1^*\cdots \tau _n^*]_F\) where \(\tau _1,\ldots
  ,\tau _n\in \Gamma ^*\) all have length \(s_0\). In fact with little
  additional effort the proof of this generalizes to our setting, and can be
  made to show that if \(s\in s_0\mathbb{N}\) then we may take all \(\tau _i\)
  to have length \(s\) rather than \(s_0\). But if we let \(a_i = [\tau _i]_F\)
  then \([\tau _1^*\cdots \tau _n^*]_F = [a_1^*\cdots a_n^*]_{F^s}\); so we
  have written \(A\) in the desired form.
\end{rpf}
We note some closure properties:
\begin{proposition}
  \label{prop:sparseclosureprops}
  \begin{menumerate}
    \item
      For any \(A\subseteq \Gamma \) and \(r>0\), \(A\) is \(F\)-sparse if and
      only if it is \(F^r\)-sparse.
    \item
      If \(A,B\subseteq \Gamma \) are \(F\)-sparse then so is \(A\cup B\).
    \item
      If \(A\subseteq \Gamma \) is \(F\)-sparse and \(X\subseteq \Gamma \) is
      \(F\)-automatic then \(A\cap X\) is \(F\)-sparse.
    \item
      If \(A,B\subseteq \Gamma \) are \(F\)-sparse then so is \(A+B\).
  \end{menumerate}
\end{proposition}
\begin{rpf}
  \begin{menumerate}
    \item
      \label{itm:pfsparseinpower}
      The right-to-left is by definition; for the left-to-right, take an
      \(F^r\)-spanning set \(\Sigma \) for some \(r>0\) and a sparse language
      \(L\subseteq \Sigma ^*\) such that \(A = [L]_{F^r}\).  Recall by
      \cite[Lemma 5.7]{bell19} that \(\Sigma ' = [\Sigma ^{(s)}]_{F^r}\) is an
      \(F^{rs}\)-spanning set.  Given \(\sigma \in\Sigma ^*\) with
      \(s\mid\abs\sigma \) we can associate \(\sigma '\in (\Sigma ')^*\) as
      follows: write \(\sigma = \sigma _1\cdots \sigma _{\frac{\abs\sigma }s}\)
      with each \(\sigma _i\in\Sigma ^{(s)}\), and set \(\sigma ' = [\sigma
      _1]_{F^r}\cdots [\sigma _{\frac{\abs\sigma }s}]_{F^r}\). So \([\sigma
      ]_{F^r} = [\sigma ']_{F^{rs}}\), and \(\abs{\sigma '} = \frac{\abs\sigma
      }s\). 

      Note now that
      \[A = [L]_{F^r}  = [L0^*\cap (\Sigma ^{(s)}))^*]_{F^r}
        = [\underbrace{\set{\sigma ' : \sigma \in L0^*\cap \Sigma
        ^{(s\mathbb{N})}}}_{L'}]_{F^{rs}}
      \]
      Since regular languages are closed under intersection (see \cite[Theorem
      2.1]{yu97}), we get that \(L0^*\cap(\Sigma ^s)^*\) is regular; one can
      then use a DFA that recognizes \(L0^*\cap(\Sigma ^s)^*\) to construct a
      DFA that recognizes \(L'\). By \cref{fact:sparsechar} \(L0^*\), and hence
      \(L0^*\cap(\Sigma ^s)^*\), is a sparse language; thus since \(\abs{\sigma
      '}= \frac{\abs\sigma }s\) we get that \(L'\) is sparse. So \(A =
      [L']_{F^{rs}}\) is \(F^s\)-sparse. 
    \item
      Take an \(F^r\)-spanning set \(\Sigma \) and an \(F^s\)-spanning set
      \(\Theta \) with sparse languages \(L_1\subseteq \Sigma ^*\) and
      \(L_2\subseteq \Theta ^*\) such that \(A = [L_1]_{F^r}\) and \(B =
      [L_2]_{F^s}\).  By the argument given in the proof of
      (\ref{itm:pfsparseinpower}), we may assume \(r=s\).  By \cite[Lemma
      5.6]{bell19}, there is an \(F^r\)-spanning set \(\Omega \) containing
      \(\Sigma \cup\Theta \). Then \(L_1\cup L_2\) is a sparse language in
      \(\Omega ^*\), and \(A\cup B = [L_1\cup L_2]_{F^r}\). So \(A\cup B\) is
      \(F\)-sparse.
    \item
      Take an \(F^r\)-spanning set \(\Sigma \) for some \(r>0\) and a sparse
      language \(L\subseteq \Sigma ^*\) such that \(A = [L]_{F^r}\). By
      \cref{prop:autoinspan} \(\set{\sigma \in\Sigma ^*:[\sigma ]_{F^r}\in X}\)
      is regular. So if \(L' = \set{\sigma \in L : [\sigma ]_{F^r}\in X}\) then
      \(L'\) is regular (as the intersection of two regular languages; see
      \cite[Theorem 2.1]{yu97}) and thus sparse (as it's contained in \(L\)).
      But \([L']_{F^r} = A\cap X\); so \(A\cap X\) is \(F\)-sparse.
    \item
      We first check the case where \(B = \set\gamma \) is a singleton. Take an
      \(F^r\)-spanning set \(\Sigma \) for some \(r>0\) and a sparse language
      \(L\subseteq \Sigma ^*\) such that \(A = [L]_{F^r}\). By \cite[Lemma
      5.6]{bell19} there is an \(F^r\)-spanning set \(\Sigma '\supseteq \Sigma
      \) that contains \(a+\gamma \) for every \(a\in\Sigma \). Given \(\sigma
      = a_1\cdots a_{\abs\sigma } \in\Sigma ^*\) non-empty let \(\sigma _\gamma
      = (a_1+\gamma )a_2\cdots a_{\abs\sigma }\in (\Sigma ')^*\); so \([\sigma
      _\gamma ]_{F^r} = [\sigma ]_{F^r}+\gamma \). It is routine to check that
      sparsity of \(L\) implies sparsity of \(L_\gamma := \set{[\sigma _\gamma
      ]_{F^r} : \sigma \in L} \subseteq (\Sigma ')^*\). So \(A+\gamma  =
      [L_\gamma ]_{F^r}\) is \(F\)-sparse.

      We now do the general case.  By \cref{prop:sparsenoconstants} there is
      \(s\) such that there is an \(F^s\)-spanning set and such that both \(A\)
      and \(B\) can be written as a finite union of translates of sets of the
      form \([a_1^*\cdots a_n^*]_{F^s}\). Distributing, it suffices to check
      the case \(A = \gamma + [a_1^*\cdots a_n^*]_{F^s}\) and \(B = \gamma
      '+[b_1^*\cdots b_{n'}^*]_{F^s}\). By the first case, we may assume
      \(\gamma =\gamma '=0\). Given \(\sigma \in a_1^*\cdots a_n^*\) and \(\tau
      \in b_1^*\cdots b_{n'}^*\) we let \(\sigma \oplus \tau \in\Gamma ^*\)
      denote the string obtained by characterwise addition; so \([\sigma \oplus
      \tau ]_{F^s} = [\sigma ]_{F^s}+[\tau ]_{F^s}\). Then
      \[a_1^*\cdots a_n^*\oplus b_1^*\cdots b_{n'}^* 
        = (a_1+b_1)^*(a_2^*\cdots a_n^*+b_1^*\cdots
        b_{n'}^*) \cup (a_1+b_1)^* (a_1^*\cdots a_n^* + b_2^*\cdots b_{n'}^*)
      \]
      and is thus sparse by an inductive argument on \((n,n')\). So \(A+B =
      [a_1^*\cdots a_n^*\oplus b_1^*\cdots b_{n'}^*]_{F^s}\) is \(F\)-sparse.
      \qedhere
  \end{menumerate}
\end{rpf}
We now work towards a characterization of \(F\)-sparsity using length
functions.
\begin{definition}
  Suppose \(A\subseteq \Gamma \) and \(\lambda \) is a length function for
  \((\Gamma ,F)\). We let \(f_{A,\lambda }(x) = \abs{\set{a\in A : \lambda
  (a)\le x}}\), and we say \(A\) is \emph{\(\lambda \)-sparse} if
  \(f_{A,\lambda }(x)\in O(\log(x)^d)\) for some \(d\in\mathbb{N}\).
\end{definition}
To explain the jump from ``polynomial in \(x\)'' to ``polynomial in
\(\log(x)\)'', recall that the length function associated to a spanning set
\(\Sigma \) was defined by raising \(2\) to the length of the shortest word
representing the input; the logarithm is there to undo the exponentiation.
Indeed, if \(\lambda =\lambda _\Sigma \) is the length function of
\cref{def:spanlenfn} associated to an \(F^r\)-spanning set \(\Sigma \) for some
\(r>0\) then \(A\) is \(\lambda \)-sparse if and only if \(\abs{\set{a\in A :
a\in [\Sigma ^{(x)}]_{F^r}}}\) grows polynomially in \(x\).

For our characterization of \(F\)-sparsity we will need the following reverse
ultrametric inequality:
\begin{lemma}
  If \(\lambda \) is a length function for \((\Gamma ,F)\) with associated
  constants \(C,D,E\) and \(a,b\in\Gamma \) satisfy \(\lambda (b) <
  D^{-1}\lambda (a)\) then \(\lambda (a+b)\ge D^{-1}\lambda (a)\).
\end{lemma}
\begin{rpf}
  Otherwise \(\lambda (a) > D\max(\lambda (b), \lambda (a+b))= D\max(\lambda
  (-b),\lambda (a+b))\ge \lambda (a)\).
\end{rpf}
We will also require the following observation relating the length of a string
to the length of the group element it represents:
\begin{lemma}
  \label{lemma:lengthboundsabstractlength}
  Suppose \(\lambda \) is a length function for \((\Gamma ,F)\) with associated
  constants \(C,D,E\); suppose \(\Sigma \) is a finite subset of \(\Gamma \).
  Then there is \(M>0\) such that \(\lambda ([\sigma ]_F)\le MC^{\abs\sigma }\)
  for all \(\sigma \in\Sigma ^*\).
\end{lemma}
\begin{rpf}
  Pick \(r\in\mathbb{N}\) such that \(C^r\ge D\), and let \(K =
  \max\set{\lambda ([\sigma ]_F) : \sigma \in\Sigma ^{(r)}}\).

  We show by induction on \(k\ge1\) that if \(\sigma \in \Sigma ^{(kr)}\) then
  \(\lambda ([\sigma ]_F)\le DC^{(k-1)r}K\). The base case is just the
  definition of \(K\). For the induction step, write \(\sigma  = \sigma
  _1\sigma _2\) where \(\abs {\sigma _2}=r\). Then by the induction hypothesis
  \(\lambda ([\sigma _1]_F)\le DC^{(k-2)r}K\), and by canonicity we get that
  \(\lambda (F^{(k-1)r}[\sigma _2]_F) \le C^{(k-1)r}\lambda ([\sigma _2]_F) \le
  C^{(k-1)r}K\).  So by the ultrametric inequality we get that \(\lambda
  ([\sigma ]_F) = \lambda ([\sigma _1]_F + F^{(k-1)r}[\sigma _2]_F) \le
  D\max(C^{(k-1)r}K, DC^{(k-2)r} K)= DC^{(k-1)r}\) (since \(C^r\ge D\)).

  Now, for \emph{any} \(\sigma \in\Sigma ^*\) we can pad by some string of
  zeroes \(\tau \) of length \(<r\) to get that \(\abs{\sigma \tau }\in
  r\mathbb{Z}\), at which point we get \(\lambda ([\sigma ]_F) = \lambda
  ([\sigma \tau ]_F) \le DC^{\abs{\sigma \tau }-r} K\le DC^{\abs \sigma }K\).
\end{rpf}
\begin{lemma}
  \label{lemma:shrinkauto}
  Suppose \(\Sigma \) is an \(F^r\)-spanning set for some \(r>0\). Suppose
  \(L\subseteq \Sigma ^*\) is regular and \(\preceq \) is a linear ordering on
  \(\Sigma \) (which induces a length-lexicographical\footnote{
    Recall that \(\sigma \preceq \tau \) in the \emph{length-lexicographical
    order} on \(\Sigma ^*\) induced by \(\preceq \) if \(\abs\sigma <\abs\tau
    \) or if \(\abs\sigma =\abs\tau \) and \(\sigma \) precedes \(\tau \) in
    the lexicographical ordering induced by \(\preceq \).
  }
  ordering on \(\Sigma ^*\) which we also denote by \(\preceq \)).  Let
  \[\widetilde{L} = \set{\sigma \in L : \sigma \preceq \tau \text{ for all
    }\tau \in L\text{ such that }[\sigma ]_{F^r} = [\tau ]_{F^r}}
  \]
  Then \(\widetilde{L}\) is regular.
\end{lemma}
\begin{rpf}
  Be replacing \(F\) with \(F^r\) we may assume \(r=1\).
  Let \(K\subseteq (\Sigma ^2)^*\) be the set of \(\col\sigma \tau \) such that
  \begin{itemize}
    \item
      \(\sigma ,\tau \in L\);
    \item
      \([\sigma ]_F = [\tau ]_F\); and
    \item
      \(\tau \prec \sigma \) or \(\tau \) ends in a \(0\).
  \end{itemize}
  So if \(\sigma \in L\) then \(\col\sigma \tau \in K\) if and only if \(\tau
  \) (with possibly the trailing \(0\) removed) witnesses that \(\sigma \notin
  \widetilde{L}\); so \(\widetilde{L} = L\setminus \pi (K)\), where \(\pi
  \colon (\Sigma ^2)^*\to\Sigma ^*\) is projection to the first coordinate. So
  since \(L\) is regular it suffices to show that \(\pi (K)\) is regular (since
  by \cite[Theorems 2.1 and 2.2]{yu97} regular languages are closed under
  Boolean combinations).

  Note that \(K\) is regular: this is because \(L\times L\), equality, ending
  in \(0\), and \(\prec \) are regular, and so \(K\) is a Boolean combination
  of regular languages. Fix a DFA \((\Sigma ^2, Q, q_0,\Omega ,\delta )\) for
  \(K\); we use this to construct an NFA \((\Sigma ,Q',q_0',\Omega ',\delta ')\)
  recognizing \(\pi (K)\). We let \(Q' = Q\), \(q_0' = q_0\), and \(\Omega ' =
  \Omega \). For the transition function we set
  \[\delta' (q, a) =\set*{\delta \pars*{q, \col ab} : b\in \Sigma }
  \]
  Then \(\sigma \) is accepted by our NFA if and only if there is \(\tau
  \in\Sigma ^*\) with \(\abs\tau =\abs\sigma \) such that \(\delta
  \pars*{q_0,\col\sigma \tau }\in\Omega \); i.e.\ such that \(\col\sigma \tau
  \in K\). So our NFA recognizes \(\pi (K)\), and \(\pi (K)\) is regular. So
  \(\widetilde{L}\) is regular.
\end{rpf}
\begin{theorem}
  \label{thm:fsparsechar}
  Suppose \(\lambda \) is a length function for some \((\Gamma ,F^r)\). Then
  \(A\subseteq \Gamma \) is \(F\)-sparse if and only if it is \(F\)-automatic
  and \(\lambda \)-sparse.
\end{theorem}
\begin{rpf}
  By \cref{prop:sparseclosureprops} (1) it suffices to check the case \(r=1\);
  so assume \(\lambda \) is a length function for \((\Gamma ,F)\). Let
  \(C,D,E\) be the constants associated to \(\lambda \).
  \begin{mdescription}
    \plr
      Suppose \(A\) is \(F\)-sparse; so there is some \(F^s\)-spanning set
      \(\Sigma \) and some sparse \(L\subseteq \Sigma ^*\) such that \(A =
      [L]_{F^s}\). We get by \cite[Proposition 6.8(b)]{bell19} that \(A\) is
      \(F\)-automatic; it remains to show that \(A\) is \(\lambda \)-sparse.
      Since \(\lambda \) is also a length function for \((\Gamma ,F^s)\) (see
      \cref{rem:lengthlargec}), we can replace \(F\) with \(F^s\) and thus
      assume that \(s=1\). By \cref{fact:sparsechar} \(L\) is a finite union of
      simple sparse languages. One can verify that a finite union of \(\lambda
      \)-sparse sets is \(\lambda \)-sparse; it thus suffices to check the case
      where \(L\) is simple sparse, say \(L = v_0w_1^*\cdots v_{n-1}w_n^*v_n\)
      with \(u_i,v_i\in\Sigma ^*\).  We apply induction on \(n\); the base case
      \(n=0\) is trivial.

      For the induction step, we have two cases.
      \begin{caselist}
        \case
          Suppose \([w_n^*v_n]_F\) is finite. Then
          \[A = \bigcup _{a \in [w_n^*v_n]_F} [v_0w_1^*\cdots
            w_{n-1}^*v_{n-1}a ]_F
          \]
          and by the induction hypothesis each \([v_0w_1^*\cdots
          w_{n-1}^*v_{n-1}a ]_F\) is \(\lambda \)-sparse. So \(A\) is \(\lambda
          \)-sparse.
        \case
          Suppose \([w_n^*v_n]_F\) is infinite.
          \begin{claim}
            There is \(M>0\) and \(i\in\mathbb{N}\) such that if \(k_n>i\) then
            \(\lambda ([v_0w_1^{k_1}\cdots w_n^{k_n}v_n]_F)\ge
            MC^{\abs{v_0w_1^{k_1}\cdots w_n^{k_n}v_n}} \).
          \end{claim}
          \begin{rpf}
            Our strategy will be to write 
            \[\underbrace{[v_0w_1^{k_1}\cdots w_n^{k_n}v_n]_F}_a =
              [v_0w_1^{k_1}\cdots w_{n-1}^{k_{n-1}}v_{n-1}w_n^{k_n-i}]_F +
              \underbrace{F^{\abs{v_0w_1^{k_1}\cdots
              w_{n-1}^{k_{n-1}}v_{n-1}w_n^{k_n-i}}}[w_n^iv_n]_F}_b
            \]
            and then use the reverse ultrametric inequality to show that 
            \(\lambda (a)\) is not much less than \(\lambda (b)\).

            Let \(M_0\) be as in \cref{lemma:lengthboundsabstractlength}. Then
            by Northcott property there is some \(i\) such that \(\lambda
            ([w_n^iv_n]_F)> EDM_0\). Suppose \(k_1,\ldots ,k_n\in\mathbb{N}\)
            with \(k_n\ge i\); to avoid notational clutter we abbreviate
            \(\sigma = v_0w_1^{k_1}\cdots
            w_{n-1}^{k_{n-1}}v_{n-1}w_n^{k_n-i}\). We then wish to show that
            \(\lambda ([\sigma w_n^iv_n]_F)\) is not much less than \(\lambda
            (F^{\abs\sigma }[w_n^iv_n]_F)\).

            Note that
            \[\lambda (F^{\abs\sigma }[w_n^iv_n]_F)\ge E^{-1}C^{\abs\sigma
              }\lambda ([w_n^iv_n]_F)> C^{\abs\sigma }DM_0\ge D \lambda
              ([\sigma ]_F)
            \]
            by hypothesis on \(M_0\).  So by the reverse ultrametric inequality
            we get that
            \[\lambda ([\sigma w_n^iv_n]_F)
              = \lambda ([\sigma ]_F + F^{\abs\sigma }[w_n^iv_n]_F)
              \ge D^{-1}\lambda (F^{\abs\sigma }[w_n^iv_n]_F)
              \ge D^{-1}E^{-1}C^{\abs\sigma }\lambda ([w_n^iv_n]_F)
            \]
            Let \(M = D^{-1}E^{-1}C^{-\abs{w_n^iv_n}}\lambda
            ([w_n^iv_n]_F)\).  Then if \(k_n\ge i\) then
            \(\lambda ([v_0w_1^{k_1}\cdots w_n^{k_n}v_n]_F)\ge
            MC^{\abs{v_0w_1^{k_1}\cdots w_n^{k_n}v_n}} \), as desired.
          \end{rpf}
          Now, we can write
          \[L = \set{v_0w_1^*\cdots w_n^{k_n}v_n : k_n \ge
            i}\cup\bigcup _{j<i} v_0w_1^*\cdots w_{n-1}^*v_{n-1}w_n^jv_n
          \]
          By the induction hypothesis each \([v_0w_1^*\cdots
          w_{n-1}^*v_{n-1}w_n^jv_n]_F\) is \(\lambda \)-sparse; it remains to
          check the case \(k_n\ge i\).

          If \(\tau =v_0w_1^{k_1}\cdots w_n^{k_n}v_n\) with \(k_n\ge i\) and
          \(\tau \) satisfies \(\lambda ([\tau ]_F)\le x\) then by the claim
          \(MC^{\abs\tau }\le x\), and thus \(\abs\tau \le \log(C)^{-1}(\log
          (x)-\log(M))\). But by hypothesis there are \(d,K\) such that
          eventually \(\abs{\set{\tau \in L :\abs\tau \le y}}\le Ky^d\). So if
          \(x\) is sufficiently large there are at most
          \(K(\log(C)^{-1}(\log(x)-\log(M))^d\in O(\log(x)^d)\) strings \(\tau
          \in L\) satisfying \(\lambda ([\tau ]_F)\le x\); so \(A\) is
          \(\lambda \)-sparse. 
      \end{caselist}
    \prl
      Let \(\Sigma \) be an \(F^s\)-spanning set for some \(s>0\). Note first
      that \(A\) is \(\lambda _\Sigma \)-sparse. Indeed, by
      \cref{lemma:lengthboundsabstractlength} (applied to \((\Gamma ,F^s)\))
      there is \(M>0\) such that \(\lambda _\Sigma ([\sigma ]_{F^s}) \le
      MC^{\abs\sigma }\) for all \(\sigma \in\Sigma ^*\).  So if \(\lambda
      _\Sigma (a)\le x\) then \(\lambda (a)\le MC^{\log_2(x)}=
      Mx^{\log_2(C)}\); thus \(f_{A,\lambda _\Sigma }(x)\le f_{A,\lambda
      }(Mx^{\log_2(C)})\in O((\log(x))^d)\) for some \(d\), and \(A\) is
      \(\lambda _\Sigma \)-sparse.

      Let \(L = \set{\sigma \in \Sigma ^* : [\sigma ]_{F^s}\in A}\); so \(L\)
      is regular by \cref{prop:autoinspan}. Fix any total order \(\preceq \) on
      \(\Sigma \), and let \(\widetilde{L}\) be as in \cref{lemma:shrinkauto};
      so \(\widetilde{L}\) is also regular. Note also that \([\cdot ]_{F^s}\)
      is a bijection \(\widetilde{L}\to A\) such that \(\lambda _\Sigma
      ([\sigma ]_{F^s}) = 2^{\abs\sigma }\): indeed, if there were a shorter
      \(\tau \in \Sigma ^*\) with \([\sigma ]_{F^s} = [\tau ]_{F^s}\) then
      \(\tau \in L\), and thus \(\tau \) witnesses that \(\sigma \notin
      \widetilde{L}\). It follows that \(\widetilde{L}\) is sparse: since \(A\)
      is \(\lambda _\Sigma \)-sparse, we get that
      \[\abs{\set{\sigma \in \widetilde{L} : \abs\sigma \le x}}
        = \abs{\set{a\in A: \lambda _\Sigma (a)\le 2^x}}
        = f_{A,\lambda_\Sigma  }(2^x)\in O(\log(2^x)^d) = O(x^d)
      \]
      for some \(d\). So \(A = [\widetilde{L}]_{F^s}\) is \(F\)-sparse.
      \qedhere
  \end{mdescription}
\end{rpf}
It follows that sparsity can be checked in any spanning set by looking at a set
of ``minimal'' representations of elements of \(A\); this can be seen as an
analogue of \cref{lemma:shrinkauto} for sparsity.
\begin{corollary}
  \label{cor:fsparsesigsparse}
  Suppose \(\Sigma \) is an \(F^r\)-spanning set for some \(r>0\). Suppose
  \(L\subseteq \Sigma ^*\) is regular and \(\preceq \) is a linear ordering on
  \(\Sigma \) (which we again identify with the induced length-lexicographic
  ordering on \(\Sigma ^*\)). Let
  \[\widetilde{L} = \set{\sigma \in L : \sigma \preceq \tau \text{ for all
    }\tau \in L\text{ such that }[\sigma ]_{F^r} = [\tau ]_{F^r}}
  \]
  Then \([L]_{F^r}\) is \(F\)-sparse if and only if \(\widetilde{L}\) is
  sparse.
\end{corollary}
\begin{rpf}
  The right-to-left direction is by definition of \(F\)-sparsity. For the
  left-to-right, note by \cref{lemma:shrinkauto} that \(\widetilde{L}\) is
  regular.  Furthermore as in the proof of \cref{thm:fsparsechar} \([\cdot
  ]_{F^r}\) is a bijection \(\widetilde{L}\to [L]_{F^r}\) with the property
  that \(\lambda _\Sigma ([\sigma ]_{F^r}) = 2^{\abs\sigma }\), and by
  \cref{thm:fsparsechar} \(A\) is \(\lambda _\Sigma \)-sparse. So as in the
  proof of \cref{thm:fsparsechar} \(\lambda _\Sigma \)-sparsity of \(A\)
  implies that \(\abs{\set{\sigma \in \widetilde{L} : \abs\sigma \le x}}\)
  grows polynomially in \(x\). So \(\widetilde{L}\) is sparse.
\end{rpf}
Another consequence is that if \(A\) is contained in an \(F\)-invariant
subgroup \(H\) then sparsity of \(A\) can be checked in \(H\):
\begin{corollary}
  \label{cor:sparsefinvsub}
  Suppose there is an \(F^r\)-spanning set for some \(r>0\). Suppose
  \(H\le\Gamma \) is \(F\)-invariant and \(A\subseteq H\). Then \(A\) is
  \(F\)-sparse in \(\Gamma \) if and only if \(A\) is \(F\)-sparse in \(H\).
\end{corollary}
\begin{rpf}
  For the right-to-left direction, note by \cref{prop:sparseclosureprops} (1)
  that \(A\) is \(F^r\)-sparse. So there is \(s>0\) with an \(F^{rs}\)-spanning
  set \(\Sigma \) for \(H\) and \(L\subseteq \Sigma ^*\) sparse such that \(A =
  [L]_{F^{rs}}\). Then using \cite[Lemmas 5.6 and 5.7]{bell19}, since there is
  an \(F^r\)-spanning set for \(\Gamma \) we may assume there is an
  \(F^{rs}\)-spanning set \(\Sigma '\) for \(\Gamma \) containing \(\Sigma \).
  Then \(L\subseteq (\Sigma ')^*\) witnesses that \(A\) is \(F\)-sparse in
  \(\Gamma \).

  For the left-to-right direction, fix any length function \(\lambda \) for
  some \((\Gamma ,F^s)\) (using \cref{thm:spanlenheight}). One can check that
  \(\lambda \) is a length function for \(H\); so by \cref{thm:fsparsechar}
  applied to \(H\) it suffices to show that \(A\) is \(F\)-automatic in \(H \)
  and \(\lambda \)-sparse. By \cref{thm:fsparsechar} applied to \(\Gamma \) we
  get that \(A\) is \(F\)-automatic in \(\Gamma \) and \(\lambda \)-sparse.
  Then \cref{cor:spanfinvsub} yields that \(A\) is \(F\)-automatic in \(H\),
  and since \(A\subseteq H\) we get that \(\lambda \)-sparsity in \(\Gamma \)
  agrees with \(\lambda \)-sparsity in \(H\). So \(A\) is \(F\)-sparse in
  \(H\).
\end{rpf}
We can now see that \(\Gamma \) (and more generally any set containing an
infinite \(F\)-invariant subgroup of \(\Gamma \)) is not \(F\)-sparse.
\begin{corollary}
  \label{cor:groupnotsparse}
  Suppose \(A\subseteq \Gamma \) is \(F\)-sparse.
  \begin{enumerate}
    \item
      \(A\) does not contain any coset of any infinite \(F\)-invariant
      subgroup.
    \item
      If \(\Gamma \) is finitely generated then \(A\) does not contain any
      coset of any infinite subgroup.
  \end{enumerate}
\end{corollary}
\begin{rpf}
  \begin{menumerate}
    \item
      We first check the case \(H=\Gamma \).  Since \(A\) is \(F\)-sparse there
      is an \(F^r\)-spanning set for some \(r>0\); by
      \cref{prop:sparseclosureprops} (1) we may assume that there is an
      \(F\)-spanning set \(\Sigma \). By \cref{lemma:canonicityinforb} we may
      assume the exceptional set of \(\lambda _\Sigma \) contains only elements
      of finite \(F\)-orbit; say the new associated constants are \(C,D,E\).
      Note that we can take \(C=2\): the constants obtained from
      \cref{prop:lenabstractlen} are \(C=D=E=2\), and
      \cref{lemma:canonicityinforb} doesn't require changing \(C\).  Note that
      there is \(a\in\Sigma \) of infinite \(F\)-orbit: otherwise by
      \cref{rem:inforderinforbit} all \(a\in\Sigma \) would be of finite order
      as well, and thus all
      \[\mathbb{Z}[F]a= \set*{\sum_{i<\abs {\set{a,Fa,\ldots }}}
        k_iF^ia : \text{ each } k_i < \abs a}
      \]
      would be finite, and \(\Gamma = \mathbb{Z}[F]\Sigma \) would be finite, a
      contradiction.

      Fix \(a\in\Sigma \) of infinite \(F\)-orbit. Fix \(s\) such that
      \(C^{s-1}>DE\).
      \begin{claim}
        Suppose \(\sigma \in\set{-a,0,a}^*\setminus \set0^* \).
        Then \([\sigma ]_{F^s}\ne 0\).
      \end{claim}
      \begin{rpf}
        We show that \(\lambda _\Sigma ([\sigma ]_{F^s})\ne 0\), which will
        suffice.

        Let \(\ell =\abs\sigma >0\). We may assume \(\sigma \) has no trailing
        zeroes and, by possibly negating, that \(\sigma _{\ell  -1}=a\); so
        \[\lambda _\Sigma ([\sigma ]_{F^s}) = \lambda_\Sigma  ([\sigma _0\cdots
          \sigma _{\ell -2}]_{F^s} + F^{s(\ell  -1)}a)
        \]
        We look to apply the reverse ultrametric inequality.  Note that
        \(\lambda _\Sigma ([\sigma _0\cdots \sigma _{\ell -2}]_{F^s})\le
        2^{s(\ell -2)+1} = C^{s(\ell -2)+1} \) by definition of \(\lambda
        _\Sigma \). Thus \(\lambda _\Sigma (F^{s(\ell -1)}a)\ge C^{s(\ell
        -1)}E^{-1}> C^{s(\ell -2)+1}D\ge D\lambda_\Sigma  ([\sigma _0\cdots
        \sigma _{\ell -2}]_{F^s})\) by canonicity and since \(C^{s-1}>DE\);
        note that canonicity applies since \(a\) has infinite \(F\)-orbit. So
        by the reverse ultrametric inequality
        \[\lambda _\Sigma ([\sigma ]_{F^s})\ge D^{-1}\lambda_\Sigma  (F^{s(\ell
          -1)}a) \ge D^{-1}E^{-1}C^{s(\ell -1)} > 0
        \]
        So in particular \([\sigma ]_{F^s}\ne0\).
      \end{rpf}
      It follows that \([\cdot ]_{F^s} \restriction \set{0,a}^*a\) is
      injective.  Indeed, suppose \(\sigma ,\tau \in\set{0,a}^*a\) are
      distinct; then we can write \([\sigma ]_{F^s}-[\tau ]_{F^s}= [\nu
      ]_{F^s}\) for some \(\nu \in\set{-a,0,a}^*\setminus \set0^*\). Hence by
      the claim \([\sigma ]_{F^s}-[\tau]_{F^s}\ne0\).

      But then 
      \[f_{\Gamma ,\lambda _\Sigma  }(2^{s\ell +1})
        = \abs{\set{[\sigma ]_F : \sigma \in\Sigma ^*, \abs\sigma \le s\ell +1}}
        \ge \abs{\set{[\sigma ]_{F^s} : \sigma \in\set{0,a}^{(\ell )}a}}
        = \abs{\set{0,a}^{(\ell )}a}
        = 2^\ell 
      \]
      So \(f_{\Gamma ,\lambda _\Sigma  }\notin O(\log(x)^d)\) for any \(d\); so
      \(\Gamma \) is not \(\lambda _\Sigma \)-sparse, and hence by
      \cref{thm:fsparsechar} \(\Gamma \) is not \(F\)-sparse.

      Suppose now that \(H\) is an arbitrary infinite \(F\)-invariant subgroup
      of \(\Gamma \). By \cref{cor:spanfinvsub} there is an \(F^s\)-spanning
      set for \(H\) for some \(s>0\); so by the case \(H=\Gamma \) we get that
      \(H\) isn't \(F\)-sparse in \(H\). So \cref{cor:sparsefinvsub} yields
      that \(H\) isn't \(F\)-sparse in \(\Gamma \).
    \item
      Again using \cref{prop:sparseclosureprops} (1) we may assume there is a
      length function \(\lambda \) for \((\Gamma ,F)\). Since \(F\)-sparsity is
      closed under translation (\cref{prop:sparseclosureprops} (2)) we may
      assume \(A\) contains an infinite \(H\le\Gamma \).  Since \(\Gamma \) is
      finitely generated and \(H\) is infinite there is some \(a\in H\) of
      infinite order; we will show that \(\mathbb{Z}a\) isn't \(\lambda
      \)-sparse. We show by induction on \(k\) that if \(1\le n\le 2^k\) then
      \(\lambda (na)\le D^{k-1}\lambda (a)\). The base case is immediate.  For
      the induction step, suppose the claim holds of \(k\); suppose \(2^k <
      n\le 2^{k+1}\). Then
      \[\lambda (na)
        = \lambda (2^ka + (n-2^k)a) \le D\max(2^ka,(n-2^k)a) \le
        D\max(D^{k-1}\lambda (a),D^{k-1}\lambda (a))\le D^k\lambda (a)
      \]
      by the induction hypothesis.

      Then \(\set{na : n\in\mathbb{Z},\lambda (na)\le D^{k-1}\lambda
      (a)}\supseteq \set{na: 1\le n\le 2^k}\) contains at least \(2^k\)
      elements. So \(f_{\mathbb{Z}a,\lambda }(x)\notin O(\log(x)^d)\) for any
      \(d\), and \(\mathbb{Z}a\) isn't \(\lambda \)-sparse. But \(A\supseteq
      H\supseteq \mathbb{Z}a\); so \(f_{A,\lambda }(x)\ge
      f_{\mathbb{Z}a,\lambda }(x)\), and \(A\) isn't \(\lambda \)-sparse. So
      \(A\) isn't \(F\)-sparse by \cref{thm:fsparsechar}.
      \qedhere
  \end{menumerate}
\end{rpf}

\section{\(F\)-sparse sets and stable expansions of finitely generated abelian
groups}
\label{sec:sparsestable}
Fix an abelian group \(\Gamma \) and a subset \(A\subseteq \Gamma \). We say
that \(A\) is \emph{stable in \(\Gamma \)} if there do not exist arbitrarily
long tuples \((a_1,\ldots ,a_N; b_1,\ldots ,b_N)\) of elements of \(\Gamma
\) such that \(a_i + b_j\in A\) if and only if \(i\le j\). This can be
expressed model-theoretically by requiring that in \(\thy(\Gamma ,0,+,A)\) the
formula \(x+y\in A\) be a stable formula. The question of which subsets of
\(\Gamma \) are stable has been of significant interest to model theorists for
some time. We refine the question here as follows:
\begin{question}
  Suppose \(\Gamma \) is a finitely generated abelian group and \(F\colon
  \Gamma \to\Gamma \) is an injective endomorphism, such that \((\Gamma ,F)\)
  admits a spanning set for some \(r>0\). Which \(F\)-automatic sets of
  \(\Gamma \) are stable?
\end{question}
In \cite{hawthorne20} I answered this question in the classical case where
\(\Gamma =\mathbb{Z}\) and \(F\) is multiplication by a positive integer. Here,
we generalize the methods of \cite{hawthorne20} to answer the above question
for \(F\)-sparse sets. This is \cref{thm:stablesparsechar} below. There are
some issues that arise in generalizing; we focus our exposition on what needs
to be done beyond what was done in \cite{hawthorne20}.

Our answer will be in terms of the ``\(F\)-cycles'' introduced in
\cite{moosa04}:
\begin{definition}
  \label{def:fsets}
  Suppose \(\Gamma \) is an abelian group and \(F\) is an injective
  endomorphism. An \emph{\(F\)-cycle} is a set of the form \(C(a;F^\delta ) :=
  \set{a + F^\delta a + \cdots + F^{\delta n}a : n<\omega }\). An
  \emph{\(F\)-set} of \(\Gamma ^m\) is a finite union of sets of the form
  \(\gamma + H+ C(a_1;F^{r_1})+\cdots + C(a_n; F^{r_n})\) where \(\gamma
  \in\Gamma ^m\), \(H\le \Gamma ^m\) is \(F\)-invariant, each \(a_i\in \Gamma
  ^m\), and each \(r_i>0\). A \emph{groupless \(F\)-set} is a finite union of
  sets of the form \(\gamma +C(a_1;F^{r_1})+\cdots + C(a_n; F^{r_n})\) with
  \(\gamma ,a_i,r_i\) as above. The \emph{\(F\)-structure} \((\Gamma
  ,\mathcal{F})\) on \(\Gamma \) has domain \(\Gamma \) and a predicate for
  every \(F\)-set of every \(\Gamma ^m\).
\end{definition}
Note in particular that the graph of addition is an \(F\)-invariant subgroup of
\(\Gamma ^3\), and is thus an \(F\)-set; so \((\Gamma ,\mathcal{F})\) expands
\((\Gamma ,+)\).

In \cite{moosa04} it was shown that if \(\Gamma \) is finitely generated and
\[\label{cond:moosascanlon}
  \bigcap _{i\in \mathbb{N}} (F^i) = \set0
  \tag{\dag}
\]
holds in \(\mathbb{Z}[F]\) then the \(F\)-structure \((\Gamma ,\mathcal{F})\)
is stable. In particular, all \(F\)-sets are stable in \(\Gamma \). Here we
prove a converse for \(F\)-sparse sets.
\begin{theorem}
  \label{thm:stablesparsechar}
  Suppose \(\Gamma \) is a finitely generated abelian group and \(F\) is an
  injective endomorphism of \(\Gamma \) such that \(\Gamma \) admits an
  \(F^r\)-spanning set for some \(r>0\).  If \(A\subseteq \Gamma \) is
  \(F\)-sparse and stable in \(\Gamma \) then \(A\) is a finite Boolean
  combination of groupless \(F\)-sets.
\end{theorem}

We would like to use the results of \cite{moosa04} to deduce quantifier
elimination of \((\Gamma ,\mathcal{F})\); unfortunately, this requires that
\(\mathbb{Z}[F]\) satisfy (\ref{cond:moosascanlon}), which doesn't follow from
the existence of a spanning set. Consider for example \(\Gamma =
\mathbb{Z}\times (\mathbb{Z}/2\mathbb{Z})\) with the endomorphism \(F(a,b) =
(2a,b)\). One can check using e.g.\ \cref{thm:spaneigen} that \(\Gamma \)
admits an \(F^r\)-spanning set for some \(r>0\), but \(F-2 = F^i(F-2)\in
(F^i)\) for all \(i\in\mathbb{N}\).

We give a sufficient condition for (\ref{cond:moosascanlon}) to hold, and then
reduce \cref{thm:stablesparsechar} to the case where this condition holds.
\begin{lemma}
  \label{lemma:moosascanloncond}
  If \(\Gamma \) is torsion-free and admits an \(F^r\)-spanning set for some
  \(r>0\) then \(\mathbb{Z}[F]\) satisfies (\ref{cond:moosascanlon}).
\end{lemma}
\begin{rpf}
  Since \(\bigcap _{i\in\mathbb{N}}(F^i) = \bigcap _{i\in\mathbb{N}}
  (F^{ri})\), we may assume there is a length function \(\lambda \) for
  \((\Gamma ,F)\) (possibly replacing \(F\) by a power thereof) with associated
  constants \(C,D,E\). By \cref{lemma:canonicityinforb,rem:inforderinforbit} we
  may assume canonicity applies to all non-zero elements.

  Suppose \(G\in\bigcap _{i\in\mathbb{N}}(F^i)\), and suppose \(a\in\Gamma \);
  suppose for contradiction that \(Ga\ne 0\). Let \(t\in\mathbb{R}\) be the
  smallest non-zero value in \(\lambda (\Gamma )\) (which must exist by the
  Northcott property). For \(n\in\mathbb{N}\) since \(G\in (F^n)\) there is
  \(G_n\in \mathbb{Z}[F]\) such that \(G = F^nG_n\); since \(Ga\ne 0\) we get
  that \(G_na\ne 0\). Then \(\lambda (Ga) = \lambda (F^nG_na) \ge
  C^nE^{-1}\lambda (G_na)\ge C^nE^{-1}t\). (Note that \(\lambda (G_na)\ne 0\)
  else canonicity would contradict Northcott.) But this grows without bound, a
  contradiction. So \(Ga = 0\) for all \(a\in\Gamma \), and \(G=0\).
\end{rpf}
We can now show \cref{thm:stablesparsechar} under the assumption that \(\Gamma
\) is torsion-free.
\begin{proposition}
  \label{prop:stablesparseinforb}
  Suppose \(\Gamma \) is torsion-free and admits an \(F^r\)-spanning set for
  some \(r>0\).  If \(A\subseteq \Gamma \) is stable and \(F\)-sparse then
  \(A\) is a finite Boolean combination of groupless \(F\)-sets.
\end{proposition}
\begin{rpf}
  By replacing \(F\) with a power thereof we may assume \(\Gamma \) admits an
  \(F\)-spanning set (since \(A\) is also \(F^r\)-sparse by
  \cref{prop:sparseclosureprops} (1)).  By \cref{prop:sparsenoconstants} we can
  write \(A\) as a finite union of sets of the form \(\gamma +[a_1^*\cdots
  a_n^*]_{F^s} \) for some \(s>0\) and \(\gamma ,a_1,\ldots ,a_n\in\Gamma \).
  If we let \(b_1,\ldots ,b_n\in \Gamma \) be such that \(a_i = \sum_{j\ge i}
  b_i\) then we can rewrite
  \[\gamma +[a_1^*\cdots a_n^*]_{F^s} = 
    \gamma + \set{[b_1^{e_1}]_{F^s} + \cdots + [b_n^{e_n}]_{F^s} : e_1\le\cdots
    \le e_n}
  \]
  and thus write \(A\) as a finite union of sets of this form.  Fix one such
  \(\gamma +\set{[b_1^{e_1}]_{F^s} + \cdots + [b_n^{e_n}]_{F^s} : e_1\le\cdots
  \le e_n}\) in the union; we show it is contained in some \(B\) of the desired
  form that is itself contained in \(A\).

  We will proceed by examining the \(e_i\) such that \([b_1^{e_1}]_{F^s} +
  \cdots + [b_n^{e_n}]_{F^s}\in A\).
  \begin{claim}
    \(X:=\set{(e_1,\ldots ,e_n)\in\mathbb{N}^n : [b_1^{e_1}]_{F^s}+\cdots +
    [b_n^{e_n}]_{F^s}\in A-\gamma }\) is quantifier-free definable in
    \((\mathbb{N},0,S,\delta \mathbb{N})\) for some \(\delta
    \in\mathbb{N}\setminus \set0\), where \(S\) is the successor function.
  \end{claim}
  \begin{rpf}
    Fix \(\sigma \in S_n\), and suppose \(e_{\sigma (1)}\le\cdots \le e_{\sigma
    (n)}\).
    Let \(c_i = \sum_{j\ge i}b_{\sigma (i)}\); so
    \[[b_1^{e_1}]_{F^s}+\cdots + [b_n^{e_n}]_{F^s}\in A\iff
      [c_1^{e_1}c_2^{e_2-e_1}\cdots c_n^{e_n-e_{n-1}}]_{F^s}\in A-\gamma 
    \]
    By \cite[Lemmas 5.6 and 5.7]{bell19} there is an \(F^s\)-spanning set
    \(\Sigma \) containing all the \(c_i\); so by \cref{prop:autoinspan}
    \(\set{\sigma \in\Sigma ^* : [\sigma ]_{F^s}\in A-\gamma }\) is regular.
    Then as argued in the proof of \cite[Proposition 2.2]{hawthorne20} there is
    some Boolean combination \(\phi (x_1,\ldots ,x_n)\) of congruences and
    equalities between one variable and one constant such that
    \[(e_1,\ldots ,e_n)\in X
      \iff [c_1^{e_1}c_2^{e_2-e_1}\cdots c_n^{e_n-e_{n-1}}]_{F^s}
      \iff (\mathbb{N},0,S,\delta \mathbb{N})\models \phi (e_1,e_2-e_1,\ldots
      ,e_n-e_{n-1})
    \]
    (again, under the assumption that \(e_{\sigma (1)}\le\cdots \le e_{\sigma
    (n)}\)). But a congruence or equality between a constant \(N\in\mathbb{N}\)
    and either \(e_1\) or \(e_{i+1}-e_i\) can be expressed as a formula in
    \((\mathbb{N},0,S,\delta \mathbb{N})\), for some sufficiently large
    \(\delta \). So this is in turn equivalent to \((\mathbb{N},0,S,\delta
    \mathbb{N})\models \psi (e_1,\ldots ,e_n)\) for some quantifier-free \(\psi
    \).

    Doing this for all \(\sigma \) and taking LCMs of the \(\delta \) and
    disjunctions over the possible orderings of the \(e_i\), we find that there
    is a single quantifier-free formula \(\phi \) in \((\mathbb{N},0,S,\delta
    \mathbb{N},<)\) such that \(X = \phi (\mathbb{N}^n)\). Then \(\phi \) is a
    stable formula under any partitioning of its variables: large ladders for
    \(\phi \) would yield large ladders for the corresponding partitioning of
    \[\bigwedge _{i=1}^n (x_i\in [b_i^*]_{F^s})\wedge (x_1+\cdots + x_n\in
      A-\gamma )
    \]
    and the latter is stable under any partitioning of the variables as \(A\)
    is stable (and addition is commutative and associative).

    But the stable quantifier-free formulas in \((\mathbb{N},0,S,\delta
    \mathbb{N},<)\) are known to be quantifier-free formulas in
    \((\mathbb{N},0,S,\delta \mathbb{N})\); this is \cite[Proposition
    3.3]{hawthorne20}. So \(X\) is quantifier-free definable in
    \((\mathbb{N},0,S,\delta \mathbb{N})\).
  \end{rpf}
  \begin{claim}
    \(Y := \set{[b_1^{e_1}]_{F^s} + \cdots + [b_n^{e_n}]_{F^s} : (e_1,\ldots
    ,e_n)\in X}\) is definable in \((\Gamma ,\mathcal{F})\).
  \end{claim}
  \begin{rpf}
    Fix \(a\in\Gamma \) such that \(F^ia\ne F^ja\) for \(i\ne j\); canonicity
    of any length function for any \((\Gamma ,F^t)\) shows such \(a\) must
    exist. Then \([a^i]_{F^s}\ne [a^j]_{F^s}\) for \(i\ne j\); indeed,
    otherwise applying \(F^s-1\) to both sides we would have \(F^{si}a=
    F^{sj}a\), a contradiction.

    Consider \(\Phi \colon \mathbb{N}\to\Gamma \) given by \(i\mapsto
    [a^i]_{F^s}\). By a similar argument to the one given in the proof of
    \cite[Lemma 3.7]{hawthorne20} we get that \(\Phi \) is an interpretation of
    \((\mathbb{N},0,S,\delta \mathbb{N})\) in \((\Gamma ,\mathcal{F})\).
    Furthermore each map \([a^i]_{F^s}\mapsto [b^i]_{F^s}\) is definable in
    \((\Gamma ,\mathcal{F})\), as is addition. So \(Y =
    \set{[b_1^{e_1}]_{F^s}+\cdots + [b_n^{e_n}]_{F^s} : ([a^{e_1}]_{F^s},\ldots
    ,[a^{e_n}]_{F^s})\in \Phi (X)}\) is definable in \((\Gamma ,\mathcal{F})\),
    as desired.
  \end{rpf}
  We get by \cref{lemma:moosascanloncond} that \((\Gamma ,F)\) satisfies
  (\ref{cond:moosascanlon}); note further since \(\Gamma \) is an infinite
  finitely generated abelian group that \(\mathbb{Z}\) embeds into
  \(\mathbb{Z}[F]\). So by \cite[Theorem A]{moosa04} \((\Gamma ,\mathcal{F})\)
  admits quantifier elimination; so \(Y\) (and hence \(\gamma +Y\)) is a
  Boolean combination of \(F\)-sets. At this point the argument given in the
  proof of \cite[Theorem 3.1]{hawthorne20} shows that \(\gamma +Y\) is a
  Boolean combination of \(F\)-sparse \(F\)-sets. Now by
  \cref{cor:groupnotsparse} an \(F\)-sparse \(F\)-set must be groupless; so
  \(\gamma +Y\) is a Boolean combination of groupless \(F\)-sets. But
  \[\gamma +\set{[b_1^{e_1}]_{F^s} + \cdots + [b_n^{e_n}]_{F^s}: e_1\le\cdots
    \le e_n}\subseteq
    \gamma +Y \subseteq A
  \]
  So we can replace \(\gamma +\set{[b_1^{e_1}]_{F^s} + \cdots +
  [b_n^{e_n}]_{F^s}: e_1\le\cdots \le e_n}\) with \(\gamma +Y\) in the union
  defining \(A\) without changing the union. Applying this to all sets in the
  union, we have written \(A\) as a Boolean combination of groupless
  \(F\)-sets.
\end{rpf}
\begin{lrpf}{thm:stablesparsechar}
  By \cref{prop:sparseclosureprops} (1) we may replace \(F\) with \(F^r\), and
  thus assume that \(\Gamma \) admits an \(F\)-spanning set \(\Sigma \).  Using
  the fundamental theorem of finitely generated abelian groups we can write
  \(\Gamma =\Gamma _0\oplus H\) where \(H\) is the torsion subgroup of \(\Gamma
  \) and \(\Gamma _0\) is torsion-free. Let \(F_0\colon \Gamma _0\to\Gamma _0\)
  be \((\pi _0\circ F)\restriction \Gamma _0\), where \(\pi _0\colon \Gamma
  \to\Gamma _0\) is the projection.  Using the fact that \(H\) is
  \(F\)-invariant one can show that \(\pi _0(\Sigma )\) is an
  \(F_0^r\)-spanning set for \(\Gamma _0\).  One can further check that if
  \(v_i,w_i\in\Sigma ^*\) then \(\pi _0([v_0w_1^{k_1}v_1\cdots
  w_n^{k_n}v_n]_{F^r}) = [\pi _0(v_0)\pi _0(w_1)^{k_1}\pi _0(v_1)\cdots \pi
  _0(w_n)^{k_n}\pi _0(w_n)]_{F_0^r}\), where for \(\sigma \in\Sigma ^*\) we
  obtain \(\pi _0(\sigma )\in(\pi _0(\Sigma))^*\) by applying \(\pi _0\) to
  each letter of \(\sigma \). So if \(B\subseteq \Gamma \) is \(F\)-sparse then
  \(\pi _0(B)\subseteq \Gamma _0\) is \(F_0\)-sparse.

  For \(b\in H\) let \(A_b = \set{a\in\Gamma _0 : a\oplus b\in A}\); note that
  \(x+y\in A_b\) is stable in \(\Gamma _0\).  We wish to show that \(A_b\)
  is \(F\)-automatic in \(\Gamma \); since regular languages are closed under
  intersection and \(A-b\) is \(F\)-automatic it suffices to show that \(\Gamma
  _0\) is.
  Fix a generating set \(\set{\gamma _1,\ldots ,\gamma _\ell }\)
  for \(\Gamma \), and write
  \[F\gamma _i = \sum_{j=1}^\ell  a_{ij}\gamma _j
  \]
  for \(a_{ij}\in\mathbb{Z}\); let \(X = (a_{ij} + \abs
  H\mathbb{Z})_{i,j=1}^\ell \in M_\ell (\mathbb{Z}/\abs H\mathbb{Z})\).  We
  construct a DFA recognizing \(\Gamma _0\). Let \(\pi _H\colon \Gamma \to H\)
  be the projection. The idea is to observe that if \(c = c_1\gamma _1+\cdots
  +c_\ell \gamma _\ell \in\Gamma \) then to compute \(\pi _H(F^nc)\) it
  suffices to know \(X^n(c_i+\abs H\mathbb{Z} )_{i=1}^\ell \), and since \(X\)
  is a matrix over a finite ring, its powers can be tracked using a finite
  automaton. Using this we can determine \(\pi _H([\sigma c]_F)\) from \(\pi
  _H([\sigma ]_F)\) and \(X^n(c_i+\abs H\mathbb{Z})_{i=1}^\ell \), which are
  finitary objects.

  The set of states is \(Q = \set{(h,Y ) : h\in H,Y\in M_\ell (\mathbb{Z}/\abs
  H\mathbb{Z})}\). The initial state is \(q_0 = (0, I)\), and the finish states
  are \(\Omega  = \set{(0,Y) : Y\in M_\ell (\mathbb{Z}/\abs H\mathbb{Z})}\).
  For the transition map, suppose we are in state \((h, Y)\) and receive input
  \(c\in \Sigma \).  Let \(Y' = XY\).  Write \(c = c_1\gamma _1+\cdots +c_\ell
  \gamma _\ell \) for \(c_i\in\mathbb{Z}\), and let \(v = (c_i + \abs
  H\mathbb{Z})_{i=1}^\ell \in (\mathbb{Z}/\abs H\mathbb{Z})^\ell \). Let
  \[h' = h + \sum_{i=1}^\ell \pi _H((Yv)_i)
  \]
  where \((Yv)_i\in \mathbb{Z}/\abs H\mathbb{Z}\) is the \(i\th\) entry of
  \(Yv\). Note that this is well-defined since \(\pi _H(\abs H\Gamma ) = 0\).
  Our machine then transitions to \((h',Y')\). One can check by induction on
  \(\abs\sigma \) that \(\delta (q_0,\sigma ) = (\pi _H([\sigma ]_F),
  X^{\abs\sigma })\).  So our machine recognizes \(\pi _H^{-1}(0)=\Gamma _0\),
  and \(\Gamma _0\) (and hence \(A_b\)) is \(F\)-automatic in \(\Gamma \).

  Furthermore \(A_b\subseteq A-b\), and \(A\) is \(F\)-sparse; so by
  \cref{prop:sparseclosureprops} (2) and (4) \(A_b\) is \(F\)-sparse in
  \(\Gamma \). So as remarked above we get that \(A_b\) is \(F_0\)-sparse in
  \(\Gamma _0\). So by \cref{prop:stablesparseinforb} each \(A_b\) is a Boolean
  combination of groupless \(F_0\)-sets.  So since
  \[A = \bigcup _{b\in H} (A_b+b)
  \]
  it suffices to show that we can write \(C(a;F_0^s)\) as a union of translates
  of sets of the form \(C(b;F^t)\). Let \(\sigma = a0^{s-1}\in\Gamma ^*\).
  Using the above automaton, we see that \(\delta (q_0,[\sigma ^i]_F)\), and
  hence \(\pi _H([\sigma ^i]_F)\), is ultimately periodic in \(i\); say \(\pi
  _H([\sigma ^{i+\mu }]_F) = \pi _H([\sigma ^i]_F)\) for \(i\ge N\). Then
  \[C(a;F_0^s) = \set{[\sigma ^i]_{F_0} : i>0}
    = \set{[\sigma ^i]_{F_0} : 0< i < N}\cup\bigcup _{j<\mu 
    }\set{[\sigma ^{N+j+i\mu }]_{F_0} : i\in\mathbb{N}}
  \]
  But \([\sigma ^{N+j+i\mu }]_F = \pi _H([\sigma ^{N+j+i\mu }]_F) + [\sigma
  ^{N+j+i\mu }]_{F_0} = \pi _H([\sigma ^{N+j}]_F) + [\sigma ^{N+j+i\mu
  }]_{F_0}\). So
  \begin{ea}
    \set{[\sigma ^{N+j+i\mu }]_{F_0} : i\in\mathbb{N}}
    &=& \pi _H([\sigma ^{N+j}]_F)+\set{[\sigma ^{N+j+i\mu }]_F :
    i\in\mathbb{N}} \\
    &=& \pi _H([\sigma ^{N+j}]_F)+[\sigma ^{N+j}]_F+
    C(F^{(N+j)s}[\sigma ^\mu ]_F;F^\mu )
  \end{ea}
  Substituting this into the above expression for \(C(a;F_0^s)\), we have
  written \(C(a;F_0^s)\) as a union of translates of \(C(b; F^\mu )\), as
  desired.
\end{lrpf}
It is worth pointing out that, assuming the existence of spanning sets, there
are always \(F\)-sparse sets that aren't \(F\)-sets. Indeed, fix any
\(a\in\Gamma \) of infinite \(F\)-orbit. One can check that \(A =
\set*{\col{F^ia}{F^ja} : i\le j}\subseteq \Gamma ^2\) is \(F\)-sparse.  If
\(A\) were an \(F\)-set then by \cite{moosa04} it would be stable,
contradicting the fact that \(\col{F^ia}0 + \col0{F^ja}\in A\) if and only if
\(i\le j\).  With some effort one can encode sets like this into subsets of
\(\Gamma \) itself.

\section{NIP expansions of \texorpdfstring{\((\Gamma ,+)\)}{(Gamma,+)}}
We now turn our attention to NIP expansions of \((\Gamma ,+)\).  We produce a
class of subsets \(A\subseteq \Gamma \), which we call the \emph{\(F\)-EDP}
sets, that contains the \(F\)-sparse sets; we show that if \(A\) is \(F\)-EDP
then \((\Gamma ,+,A)\) is NIP. Our argument generalizes that of
\cite[Subsection 6.2]{hawthorne20}, which deals with the case when
\(\Gamma =\mathbb{Z}\) and \(F\) is multiplication by some \(d>1\).  As an
application of our general result, we will see that
\((\mathbb{F}_p[t],+,\tpowmult)\) is NIP, where \(\tpowmult\) is the
graph of multiplication on \(t^\mathbb{N}\).

We introduce some convenient multi-index notation. Suppose \(\mathbf{s} =
(\sigma _1,\ldots ,\sigma _n)\) with \(\sigma _1,\ldots ,\sigma _n\in\Gamma
^*\) and \(\mathbf{k} = (k _1,\ldots ,k _n)\in\mathbb{N}^n\). Then by
\(\mathbf{s}^{\mathbf{k}}\) we mean the string \(\sigma _1^{k_1}\cdots
\sigma _n^{k_n}\). Note that if \(\mathbf{0} = (0,\ldots ,0)\) then
\(\mathbf{s}^{\mathbf0} = \epsilon \), the empty word.
\begin{definition}
  By an \emph{\(F\)-exponentially definable in Presburger} (\(F\)-EDP) subset
  of \(\Gamma \) we mean a set of the form \([\mathbf{s}^{\phi(\mathbb{N})} ]_F
  := \set*{ [\mathbf{s}^{\mathbf{k}}]_F : (\mathbb{N},+)\models \phi
  (\mathbf{k})} \) for some tuple \(\mathbf{s}\) of strings over \(\Gamma \)
  and some formula \(\phi (\mathbf{x})\) with \(\abs{\mathbf{x}} =
  \abs{\mathbf{s}}\).
\end{definition}
\begin{lemma}
  \label{lemma:edpprops}
  \begin{menumerate}
    \item
      \(F\)-EDP sets are closed under finite union.
    \item
      \(F\)-sparse sets are \(F\)-EDP.
  \end{menumerate}
\end{lemma}
\begin{rpf}
  \begin{menumerate}
    \item
      It suffices to check pairwise unions. Suppose we are given \(F\)-EDP sets
      \([\mathbf{s}_1^{\phi (\mathbb{N})}]_F\) and \([\mathbf{s}_2^{\psi
      (\mathbb{N})}]_F\).  Then their union can be written as
      \([(\mathbf{s}_1\mathbf{s}_2)^{\chi (\mathbb{N})}]_F\) where \(\chi
      (\mathbf{x},\mathbf{y})\) is \( (\phi (\mathbf{x})\wedge
      \mathbf{y}=\mathbf{0})\vee (\psi (\mathbf{y})\wedge
      \mathbf{x}=\mathbf{0})\), and is thus \(F\)-EDP.
    \item
      Suppose we are given an \(F^r\)-spanning set \(\Sigma \) for some \(r>0\)
      and a sparse \(L\subseteq \Sigma ^*\).  By \cref{fact:sparsechar} we may
      assume \(L\) is a finite union of simple sparse languages; by part (1) we
      may assume \(L\) itself is simple sparse, say \(L = v_0w_1^*v_1\cdots
      w_n^*v_n\) for \(v_i,w_i\in\Sigma ^*\).  Given \(\sigma = s_0\cdots
      s_{\abs\sigma-1}\in \Sigma ^*\) let \(\sigma ' =
      s_00^{r-1}s_10^{r-1}\cdots s_{\abs\sigma -1}0^{r-1}\).  Then \((\sigma
      \tau)' = \sigma '\tau '\) and \([\sigma ']_F = [\sigma ]_{F^r}\). So
      \([L]_{F^r} = \set{[v_0'(w_1')^{k_1}v_1'\cdots (w_n')^{k_n}v_n']_F :
      k_1,\ldots ,k_n\in\mathbb{N}}\) is \(F\)-EDP.
      \qedhere
  \end{menumerate}
\end{rpf}
The \(F\)-EDP sets contain significantly more than just the \(F\)-sparse sets.
For instance, if \(a,b,c\in\Gamma \) then \(\set{[a^ib^i]_F : i
\in\mathbb{N}}\) and \(\set{[a^ib^jc^{i+j}]_F : i,j\in\mathbb{N}}\) are
\(F\)-EDP, but are not typically \(F\)-automatic, and hence not \(F\)-sparse.

Our goal is to prove that if we start with a weakly minimal abelian group
\((\Gamma,+) \), and if \(A\subseteq \Gamma \) is \(F\)-EDP, then \((\Gamma
,+,A)\) is an NIP structure. Recall that \(\Gamma \) being weakly minimal means
that \(\thy(\Gamma ,+)\) is superstable and of \(U\)-rank \(1\). Equivalently,
for all \(n>0\), \(n\Gamma \) and the subgroup of \(n\)-torsion are either
finite or of finite index in \(\Gamma \)---see for example \cite[Proposition
3.1]{conant20}. So, for example, all finitely generated abelian groups are
weakly minimal.

We begin by adapting \cref{prop:sparsenoconstants} to EDP sets.
\begin{lemma}
  \label{lemma:edpform}
  If \(A\) is \(F\)-EDP then there is \(s_0\) such that for all \(s\in
  s_0\mathbb{N}\) we can write \(A\) in the form \([\mathbf{a}^{\phi
  (\mathbb{N})}]_{F^s}\) for some formula \(\phi \) in Presburger arithmetic
  and some tuple \(\mathbf{a}\) of elements of \(\Gamma \) (i.e.\ strings of
  length \(1\)).
\end{lemma}
\begin{rpf}
  Write \(A = \set{[\mathbf{s}^{\mathbf{k}}]_F:(\mathbb{N},+)\models \phi
  (\mathbf{k})}\) where \(\mathbf{s}= (\sigma _1,\ldots ,\sigma _n)\). Let
  \(s_0 = \lcm(\abs{\sigma _1},\ldots ,\abs{\sigma _n})\), and suppose \(s\in
  s_0\mathbb{N}\).  If \(\sigma _i' = \sigma _i^{\frac s{\abs{\sigma
  _i}}}\) then
  \[A = \bigcup _{i_1 < \frac s{\abs{\sigma _1}}}\cdots \bigcup _{i_n < \frac
    s{\abs{\sigma _n}}} \set*{[(\sigma _1')^{k_1'}\sigma _1^{i_1}\cdots (\sigma
      _n') ^{k_n'}\sigma _n ^{i_n }]_F : (\mathbb{N},+)\models \phi
      \pars*{\frac s{\abs{\sigma _1}} k_1'+i_1,\ldots ,\frac s{\abs {\sigma _n
    }}k_n '+i_n }}
  \]
  Since sets of the desired form are closed under union (as in
  \cref{lemma:edpprops} (1)) it thus suffices to check
  the case
  \[A = \set{[v _0w_1^{k_1}\cdots w _n ^{k_n }v _n]_F
    : (\mathbb{N},+)\models \phi (k_1,\ldots ,k_n )}
  \]
  where the \(w _i\) all have the same length \(s\). But by the argument given
  in \cite[Lemma 3.5]{hawthorne20} (which generalizes to our context, as we
  noted in the proof of \cref{prop:sparsenoconstants}) there are \(\gamma
  \in\Gamma \) and \(\tau _1,\ldots ,\tau _n \in \Gamma^* \) of length \(s\)
  such that if each \(k_j>0\) then we can write
  \[[v_0w _1^{k_1}\cdots w _n ^{k_n }v _n ]_F =
    \gamma + [\tau _1^{k_1-1}\cdots \tau _n ^{k_n -1}]_F
    =\gamma + [a_1^{k_1-1}\cdots a_n^{k_n-1}]_{F^s}
  \]
  where \(a_i = [\tau _i]_F\).
  The case where some \(k_j=0\) can be dealt with inductively, again using
  closure under union; so it suffices to check the case \(\gamma +
  \set{[a_1^{k_1}\cdots a_n^{k_n}]_{F^s} : (\mathbb{N},+)\models \phi
  (k_1,\ldots ,k_n)}\). But we can rewrite this as
  \[\bigcup _{i=1}^n \set{[(\gamma +a_i)a_i^{k_i}\cdots a_n^{k_n}]_{F^s} :
    (\mathbb{N},+)\models \phi (0,\ldots ,0,k_i+1,k_{i+1},\ldots ,k_n)}
  \]
  which takes the desired form by closure under union.
\end{rpf}
The source of NIP in \cref{thm:edp1dim} will be that \((\mathbb{N},+)\) is NIP.
The following lemma and corollary relate automatic and \(F\)-EDP sets to sets
definable in \((\mathbb{N},+)\).
\begin{lemma}
  Suppose \(\Lambda \) is a finite alphabet and \(L\subseteq (\Lambda^m) ^*\)
  is regular. Suppose \( \mathbf{a}_1,\ldots ,\mathbf{a}_m\) are tuples from
  \(\Lambda \).  Then the relation
  \[\set*{(\mathbf{k}_1,\ldots ,\mathbf{k}_m) : 
      \abs{\mathbf{a}_1^{\mathbf{k}_1}}= \cdots =
      \abs{\mathbf{a}_m^{\mathbf{k}_m}},
      \tcol{\mathbf{a}_1^{\mathbf{k}_1}}\vdots {\mathbf{a}_m^{\mathbf{k}_m}}\in
    L}
    \subseteq \mathbb{N}^{\abs{\mathbf{a}_1}}\times \cdots
    \times \mathbb{N}^{\abs{\mathbf{a}_m}}
  \]
  is definable in \((\mathbb{N},+)\). (Here we identify \((\Lambda ^m)^*\) with
  the subset of \((\Lambda ^*)^m\) of tuples whose constituent strings all have
  the same length.)
\end{lemma}
\begin{rpf}
  Fix an automaton \((\Lambda ^m,Q,q_0,\Omega,\delta  )\) for \(L \); we show
  that for any \(q_1,q_2\in Q\) the relation
  \[\set*{(\mathbf{k}_1,\ldots ,\mathbf{k}_m) : 
      \abs{\mathbf{a}_1^{\mathbf{k}_1}}= \cdots =
      \abs{\mathbf{a}_m^{\mathbf{k}_m}},
      \delta \pars*{q_1,\tcol{\mathbf{a}_1^{\mathbf{k}_1}}\vdots
    {\mathbf{a}_m^{\mathbf{k}_m}}} = q_2}
  \]
  is definable in \((\mathbb{N},+)\). Applying this to \(q_0\) and the elements
  of \(\Omega \) yields the desired result.

  We apply induction on \(\abs{\mathbf{a}_1}\cdots \abs{\mathbf{a}_m}\). For
  the base case, if some \(\abs{\mathbf{a}_i} = 0\) then our relation is either
  empty or equivalent to all \(\mathbf{k}_i\) being zero, depending on whether
  \(q_1 = q_2\), and both of these are definable in \((\mathbb{N},+)\).

  For the induction step, suppose no \(\abs{\mathbf{a}_i} = 0\). Write
  \(\mathbf{a}_i = (a_{i1},\ldots ,a_{in_i})\) and \(\mathbf{k}_i =
  (k_{i1},\ldots ,k_{in_i})\). Without loss of generality suppose \(k_{11}\) is
  minimum among the \(k_{i1}\). If
  \[q = \delta \pars*{q_1, \tcol{a _{11}^{k_{11}}}{\vdots }{a
    _{m1}^{k_{11}}}}
  \]
  then our relation is equivalent to
  \[
        \abs{ a _{12}^{k_{12}}\cdots a _{1n_1}^{k_{1n_1}}}
        = \abs{
        a _{21}^{k_{21}-k_{11}}a _{22}^{k_{22}}\cdots
        a _{2n_2}^{k_{2n_2}}
      }
      =\cdots =\abs{
        a _{m1}^{k_{m1}-k_{11}}a _{m2}^{k_{m2}}\cdots
        a _{mn_m}^{k_{mn_m}}
      }
    \wedge 
    \delta \pars*{q, \fcol{
      \hphantom{a _{11}^{k_{11}-k_{11}}}a _{12}^{k_{12}}\cdots a
      _{1n_1}^{k_{1n_1}}
    }{
      a _{21}^{k_{21}-k_{11}}a _{22}^{k_{22}}\cdots
      a _{2n_2}^{k_{2n_2}}
    }\vdots {
      a _{m1}^{k_{m1}-k_{11}}a _{m2}^{k_{m2}}\cdots
      a _{mn_m}^{k_{mn_m}}
    }} = q_2
  \]
  which by the induction hypothesis is definable in \((\mathbb{N},+)\).
  Similarly we get definability in the case \(k_{i1}\) is minimum for some
  \(i>1\). So taking disjunctions we get that our relation is definable in
  \((\mathbb{N},+)\).
\end{rpf}
\begin{corollary}
  \label{cor:automatonpresdfbl}
  Suppose \(\Gamma \) admits an \(F\)-spanning set. If \(X\subseteq \Gamma
  ^m\) is \(F\)-automatic and \(\mathbf{a}_1,\ldots
  ,\mathbf{a}_m\) are tuples from \(\Gamma \) then
  \[\set*{(\mathbf{k}_1,\ldots ,\mathbf{k}_m) : 
    \tcol{\relax  [\mathbf{a}_1^{\mathbf{k}_1}]_F}\vdots {\relax
    [\mathbf{a}_m^{\mathbf{k}_m}]_F} \in X}\subseteq
    \mathbb{N}^{\abs{\mathbf{a}_1}}\times \cdots
    \times \mathbb{N}^{\abs{\mathbf{a}_m}}
  \]
  is definable in \((\mathbb{N},+)\).
\end{corollary}
\begin{rpf}
  By \cite[Lemma 5.6]{bell19} there is an \(F\)-spanning set \(\Sigma \)
  containing all the \(a_{ij}\); one can then check that \(\Sigma ^m\) is an
  \(F\)-spanning set for \(\Gamma ^m\).  Let \(L\subseteq (\Sigma ^m)^*\) be
  the set of representations of elements of \(X\); so \(L\) is regular. Then by
  previous lemma
  \[\set*{(\mathbf{k}_1,\ldots ,\mathbf{k}_m,\ell _1,\ldots ,\ell _m) : 
      \abs{\mathbf{a}_1^{\mathbf{k}_1}0^{\ell _1}} = \cdots =
      \abs{\mathbf{a}_m^{\mathbf{k}_m}0^{\ell _m}},
      \tcol{\mathbf{a}_1^{\mathbf{k}_1}0^{\ell _1}}\vdots {
    \mathbf{a}_m^{\mathbf{k}_m}0^{\ell _m}} \in L}
  \]
  is definable in \((\mathbb{N},+)\), say by \(\phi (\mathbf{k}_1,\ldots
  ,\mathbf{k}_m,\ell _1,\ldots ,\ell _m)\). Then \(\exists \ell _1\cdots
  \exists \ell _m\phi (\mathbf{k}_1,\ldots ,\mathbf{k}_m,\ell _1,\ldots ,\ell
  _m)\) defines the desired relation in \((\mathbb{N},+)\).
\end{rpf}
In order to apply the previous two results, we will need to relate the
relations we wish to show are NIP to automatic sets. The following lemma does
so. It is a straightforward generalization of a known result in classical
automata theory; see for example \cite[Theorem 6.1]{bruyere94}.
\begin{lemma}
  \label{lemma:dfblautomatic}
  Suppose \(\Gamma \) admits an \(F^r\)-spanning set for some \(r>0\).
  If \(X\subseteq \Gamma ^m\) is definable in \((\Gamma ,+)\) with parameters
  from \(\Gamma \) then \(X\) is \(F\)-automatic.
\end{lemma}
\begin{rpf}
  We apply structural induction on formulas. That the claim holds for \(x = y\)
  and \(x+y = z\) is \cite[Lemma 6.7]{bell19}; that the claim holds for a union
  or complement is because regular languages are closed under union and
  complementation. It remains to check existential quantification.

  Suppose then that \(X\subseteq \Gamma ^{m+1}\) is \(F\)-automatic. Fix an
  \(F^r\)-spanning set \(\Sigma \) for some \(r>0\); then as before one can
  check that \(\Sigma ^{m+1}\) is an \(F^r\)-spanning set for \(\Gamma
  ^{m+1}\).  So by \cref{prop:autoinspan} there is an automaton \(M = (\Sigma
  ^{m+1},Q,q_0,\Omega ,\delta )\) recognizing the representations over \(\Sigma
  ^{m+1}\) of elements of \(X\).  Consider the non-deterministic automaton \(M'
  = (\Sigma ^m, Q,q_0,\Omega ',\delta ')\) where
  \begin{itemize}
    \item
      \(\delta '(q,\mathbf{a}) = \set*{\delta \pars*{q,\col {\mathbf{a}}b}
      : b\in\Sigma }\) for \(\mathbf{a}\in\Sigma ^m\), and
    \item
      \(\Omega ' = \set*{q\in Q : \text{there is }\tau \in\Sigma ^*\text{ such
      that }\delta \pars*{q,\col{\mathbf{0}^{\abs\tau }}\tau }\in \Omega }\).
  \end{itemize}
  Then if \(\boldsymbol{\sigma }\in(\Sigma ^m)^*\) is accepted by \(M'\) then
  there is \(\tau_1 \in\Sigma ^{(\abs{\boldsymbol{\sigma}}) }\) such that
  \(\delta \pars*{q_0,\col{\boldsymbol{\sigma }}{\tau _1}}\in \Omega '\); i.e.\
  such that there is \(\tau _2\in\Sigma ^*\) such that \(\delta
  \pars*{q_0,\col{\boldsymbol{\sigma }\mathbf{0}^{\abs{\tau _2}}}{\tau _1\tau
  _2}}\in \Omega \).  So \(\col{\relax [\boldsymbol{\sigma} ]_{F^r}}{\relax
  [\tau _1\tau _2]_{F^r}}\in X\). Conversely suppose there is \(b\in\Gamma \)
  such that \(\col{\relax [\boldsymbol{\sigma }]_{F^r}}b\in X\); say \(b =
  [\tau ] _{F^r}\) for \(\tau \in\Sigma ^*\). Write \(\tau =\tau _1\tau _2\)
  where \(\abs{\tau _1} = \abs{\boldsymbol{\sigma }}\) (replacing \(\tau \)
  with some \(\tau 0^k\) if necessary). Then \(\col{\relax [\boldsymbol{\sigma
  }\mathbf{0}^{\abs{\tau _2}}]_{F^r}}{\relax [\tau _1\tau _2]_{F^r}} = \col
  {\relax [\boldsymbol{\sigma} ]_{F^r}}b\in X\); so \(\col{\boldsymbol{\sigma
  }\mathbf{0}^{\abs{\tau _2}}}{\tau _1\tau _2}\) is accepted by \(M\).  So
  \(\delta \pars*{q_0,\col{\boldsymbol{\sigma }}{\tau _1}} \in \delta
  '(q_0,\boldsymbol{\sigma })\cap \Omega '\), and \(M'\) accepts
  \(\boldsymbol{\sigma }\).

  So \(M'\) recognizes the representations over \(\Sigma ^r\) of the projection
  of \(X\) away from the last coordinate.  So existentially quantifying an
  \(F\)-automatic set results in an \(F\)-automatic set.
\end{rpf}
\begin{theorem}
  \label{thm:edp1dim}
  Suppose \((\Gamma,+) \) is a weakly minimal abelian group and admits an
  \(F^r\)-spanning set for some \(r>0\); suppose \(A\subseteq \Gamma \) is
  \(F\)-EDP. Then \((\Gamma ,+,A)\) is NIP.
\end{theorem}
\begin{rpf}
  We apply \cite[Theorem 2.9]{conant20}, which tells us that in a weakly
  minimal group \(\Gamma \) if \(A_\Gamma \) is NIP then so is \((\Gamma
  ,+,A)\). Here \(A_\Gamma \) is the \emph{induced structure} of \((\Gamma,+)
  \) on \(A\): the domain of \(A_\Gamma \) is \(A\), and for every \((\Gamma
  ,+)\)-definable subset \(X\) of \(\Gamma ^m\) with parameters from \(\Gamma
  \) its trace \(X\cap A^m\) on \(A\) is an atomic relation of \(A_\Gamma \).
  We show that \(A_\Gamma \) is interpretible in \((\mathbb{N},+)\), and thus
  is NIP.

  By \cref{lemma:edpform} we may assume \(A\) takes the form
  \([\mathbf{a}^{\phi (\mathbb{N})}]_{F^s}\) for some \(s\in r\mathbb{N}\). So
  \(\Gamma \) admits an \(F^s\)-spanning set by \cite[Lemma 5.7]{bell19}. Let
  \(\Phi \colon \mathbb{N}^{\abs{\mathbf{a}}}\to \Gamma \) be
  \(\mathbf{k}\mapsto [\mathbf{a}^{\mathbf{k}}]_{F^s} \);  so \(\Phi \) is
  surjective \(\phi (\mathbb{N}^{\abs{\mathbf{a} }})\twoheadrightarrow A\). We
  show that \(\Phi \) defines an interpretation of \(A_\Gamma \) in
  \((\mathbb{N},+)\), and hence that \(A_\Gamma \) is NIP.

  By \cref{cor:automatonpresdfbl} the equivalence relation given by
  \(\mathbf{k}_1\sim \mathbf{k}_2\) if and only if \(\Phi (\mathbf{k}_1) = \Phi
  (\mathbf{k}_2)\) is definable in \((\mathbb{N},+)\) (since the diagonal in
  \(\Gamma ^{2m}\) is \(F\)-automatic by \cite[Lemma 6.7]{bell19}). Suppose now
  that \(X\subseteq \Gamma ^m\) is definable with parameters in \((\Gamma
  ,0,+)\); it remains to show that \( \phi
  (\mathbb{N}^{\abs{\mathbf{a}}})^m\cap (\Phi ^m)^{-1}(X) \)
  is definable in \((\mathbb{N},+)\).  We get by \cref{lemma:dfblautomatic}
  that \(X\) is \(F\)-automatic; so by \Cref{cor:automatonpresdfbl} \((\Phi
  ^m)^{-1}(X)\) is definable in \((\mathbb{N},+)\). So \(\phi
  (\mathbb{N}^{\abs{\mathbf{a}}})^m\cap (\Phi ^m)^{-1}(X)\) is definable in
  \((\mathbb{N},+)\), and \(\Phi \) defines an interpretation of \(A_\Gamma \)
  in \((\mathbb{N},+)\).
\end{rpf}
As an application of \cref{thm:edp1dim}, we prove the following:
\begin{theorem}
  \label{thm:polysnip}
  Suppose \(p\ge 7\) is prime. Then \((\mathbb{F}_p[t],+,\tpowmult)\) is NIP,
  where \(\tpowmult = \set*{\tcol{t^i}{t^j}{t^{i+j}} : i,j\in\mathbb{N}}\).
\end{theorem}
\begin{rpf}
  Let \(F\colon \mathbb{F}_p[t]\to \mathbb{F}_p[t]\) be \(f(t)\mapsto tf(t)\).
  Note that \(\mathbb{F}_p\) is an \(F\)-spanning set for \(\mathbb{F}_p[t]\).
  Indeed, an arbitrary element of \(\mathbb{F}_p[t]\) is of the form \(a_0 +
  a_1t+\cdots + a_nt^n = [a_0a_1\cdots a_n]_F\), and the other axioms are
  easily verified using the facts that \(\mathbb{F}_p\) is a subgroup and
  \(\mathbb{F}_p\cap F(\mathbb{F}_p[t]) = \set0\).  Furthermore
  \((\mathbb{F}_p[t],+)\) is weakly minimal: for each \(n\in\mathbb{N}\) we
  have
  \begin{ea}
    n\mathbb{F}_p[t] &=& 
    \begin{cases}
      \mathbb{F}_p[t] &\tif p\nmid n \\
      0 &\telse
    \end{cases} \\
    \set{a\in \mathbb{F}_p[t] : na = 0} &=&
    \begin{cases}
      0 &\tif p\nmid n \\
      \mathbb{F}_p[t] &\telse
    \end{cases}
  \end{ea}
  so weak minimality follows from \cite[Proposition 3.1]{conant20}. So
  \cref{thm:edp1dim} applies, and \((\mathbb{F}_p[t],+,A)\) is NIP whenever
  \(A\subseteq \mathbb{F}_p[t]\) is \(F\)-EDP.

  We encode \(\tpowmult\) in an \(F\)-EDP subset of \(\mathbb{F}_p[t] \). Write
  \(\mathbb{F}_p = \mathbb{Z}/p\mathbb{Z} = \set{\overline{n}:
  n\in\mathbb{Z}}\).  Let
  \[A = t^\mathbb{N} \cup \overline{2}t^\mathbb{N} \cup
    \set{\overline{3}t^{i+j} - \overline{3}t^i - \overline{3}t^j :
    i,j\in\mathbb{N}\setminus \set0}
  \]
  Note that
  \begin{ea}
    \set{\overline{3}t^{i+j}-\overline{3}t^i-\overline{3}t^j : i,j\in
    \mathbb{N}\setminus \set0}
    &=& \set{[(\overline{0})^k(\overline{-1})(\overline{0})^\ell
    (\overline{-1})(\overline{0})^{k-1}\overline{1}]_F : k,\ell
    \in\mathbb{N},k>0}\\
    &&\cup
    \set{[(\overline{0})^k(\overline{-2})(\overline{0})^{k-1}\overline{1}]_F :
    k\in\mathbb{N}\setminus \set0}
  \end{ea}
  is \(F\)-EDP as a union of \(F\)-EDP sets (see \cref{lemma:edpprops} (1)).
  Since \(t^\mathbb{N} = \set{[(\overline{0})^k\overline{1}]_F : k
  \in\mathbb{N}}\), we get that \(A\) itself is \(F\)-EDP as a union of
  \(F\)-EDP sets; so \((\mathbb{F}_p[t] ,+,A)\) is NIP. It remains to see that
  \(\tpowmult\) is definable in \((\mathbb{F}_p[t],+,A)\).
  \begin{claim}
    \(t^\mathbb{N}\) and \(B := \set{t^{i+j} - t^i - t^j :
    i,j\in\mathbb{N}\setminus \set0}\) are both definable in \((\mathbb{F}_p[t]
    ,+,A)\).
  \end{claim}
  \begin{rpf}
    We first check \(t^\mathbb{N}\). I claim that if \(\phi (x)\) is \( (x\in
    A)\wedge (2x\in A)\) then \(\phi \) defines \(t^\mathbb{N}\). It is clear
    that \(\phi (\mathbb{F}_p[t] )\supseteq t^\mathbb{N}\); suppose conversely
    that \(f\in\phi (\mathbb{F}_p[t] )\). If we had \(f = \overline{2}t^i\in
    \overline{2}t^\mathbb{N}\) then \(2f = \overline{4}t^i\), and hence
    \(2f\notin A\) (since \(p\ge7\) implies \(\overline{4}\notin
    \set{\overline{1},\overline{2},\overline{3}}\), and these are the only
    leading coefficients of elements of \(A\)). But this contradicts the
    assumption that \(f\in\phi (\mathbb{F}_p[t])\).  Similarly if we had \(f=
    \overline{3}t^{i+j}-\overline{3}t^i-\overline{3}t^j\) for some \(i,j\) then
    \(2f = \overline{6}t^{i+j} - \overline{6}t^i - \overline{6}t^j\) has
    leading coefficient \(\overline{6}\); so by similar reasoning we conclude
    that \(2f\notin A\), a contradiction.  So having eliminated the other
    possibilites we conclude that \(f\in t^\mathbb{N}\), as desired. Note then
    that \(\overline{2}t^\mathbb{N} = 2t^\mathbb{N}\) is also definable.

    We now check \(B\). I claim that if \(\psi (x)\) is \(3x\in A\setminus
    (t^\mathbb{N}\cup 2t^\mathbb{N})\) then \(\psi \) defines \(B\). That
    \(\psi (\mathbb{F}_p[t] )\subseteq B\) is just injectivity of \(f\mapsto
    3f\); that \(\psi (\mathbb{F}_p[t] )\supseteq B\) is that
    \(\set{\overline{3}t^{i+j}-\overline{3}t^i-\overline{3}t^j :
    i,j\in\mathbb{N}}\cap (t^\mathbb{N}\cup\overline{2}t^\mathbb{N}) =
    \emptyset \), which follows by considering the leading coefficients.
  \end{rpf}
  So it suffices to show that \(\tpowmult\) is definable in \((\mathbb{F}_p[t]
  ,+,t^\mathbb{N},B)\).  In fact if \(\phi (x,y,z)\) is
  \[(x = \overline{1}\wedge z = y\in t^\mathbb{N})\vee (y=\overline{1} \wedge
    z=x\in t^\mathbb{N})\vee (x,y,z\in t^\mathbb{N}\wedge z-x-y \in B )
  \]
  then \(\tpowmult = \phi (\mathbb{F}_p[t] )\). It is clear that
  \(\tpowmult\subseteq \phi (\mathbb{F}_p[t] )\); for the converse suppose
  \(\tcol{t^i}{t^j}{t^k}\in\phi (\mathbb{F}_p[t] )\). Note that if either of
  the first two disjuncts of \(\phi \) holds then \(\tcol{t^i}{t^j}{t^k}\in
  \tpowmult\); suppose then that the last holds. So there are \(i',j'>0\) such
  that \(t^k-t^i-t^j = t^{i'+j'}-t^{i'} -t^{j'} \). Since \(i',j'>0\) we get
  that the leading coefficient of the right-hand side is \(1\); for the same to
  be true on the left-hand side we must have \(k >i,j\), and thus \(k =
  i'+j'\). So \(t^i + t^j = t^{i'} + t^{j'}\), and \(\set{i,j} = \set{i',j'}\).
  Thus \(i+j = i'+j' = k\), and \(\tcol{t^i}{t^j}{t^k}\in \tpowmult\).
\end{rpf}
\printbibliography
\end{document}